\documentclass{article}
\usepackage{amsmath}
\usepackage{amssymb}

\usepackage{theorem}
\numberwithin{equation}{section}

\newtheorem{theorem}{Theorem}[section]

\newtheorem{lemma}[theorem]{Lemma}

\begin{document}
\title{Global well-posedness for the radial, defocusing, nonlinear wave equation for $3 < p < 5$}
\date{\today}
\author{Benjamin Dodson}
\maketitle

\noindent \textbf{Abstract:} In this paper we continue the study of the defocusing, energy-subcritical nonlinear wave equation with radial initial data lying in the critical Sobolev space. In this case we prove scattering in the critical norm when $3 < p < 5$.

\section{Introduction}
In this paper we study the defocusing, nonlinear wave equation
\begin{equation}\label{1.1}
u_{tt} - \Delta u + |u|^{p - 1} u = 0, \qquad u(0,x) = u_{0}, \qquad u_{t}(0,x) = u_{1},
\end{equation}
for $3 < p < 5$. This problem has the critical scaling symmetry
\begin{equation}\label{1.2}
u(t,x) \mapsto \lambda^{\frac{2}{p - 1}} u(\lambda t, \lambda x).
\end{equation}
Under this scaling, the critical Sobolev exponent
\begin{equation}\label{1.3}
s_{c} = \frac{3}{2} - \frac{2}{p - 1}
\end{equation}
is preserved. In this paper we continue the study that we began in \cite{D}, proving
\begin{theorem}\label{t1.1}
The defocusing, nonlinear wave equation
\begin{equation}\label{1.4}
u_{tt} - \Delta u + |u|^{p - 1} u = 0, \qquad u(0,x) = u_{0}, \qquad u_{t}(0,x) = u_{1},
\end{equation}
is globally well-posed and scattering for radial initial data $(u_{0}, u_{1}) \in \dot{H}^{s_{c}}(\mathbf{R}^{3}) \times \dot{H}^{s_{c} - 1}(\mathbf{R}^{3})$. Moreover, there exists a function $f : [0, \infty) \rightarrow [0, \infty)$ such that if $u$ solves $(\ref{1.4})$ with initial data $(u_{0}, u_{1}) \in \dot{H}^{s_{c}} \times \dot{H}^{s_{c} - 1}$, then
\begin{equation}\label{1.4.1}
\| u \|_{L_{t,x}^{2(p - 1)}(\mathbf{R} \times \mathbf{R}^{3})} \leq f(\| u_{0} \|_{\dot{H}^{s_{c}}} + \| u_{1} \|_{\dot{H}^{s_{c} - 1}}).
\end{equation}
\end{theorem}

There are a number of reasons to conjecture that such a result is true. First, it is known that critical Sobolev regularity completely determines local well-posedness.
\begin{theorem}\label{t1.2}
The equation $(\ref{1.4})$ is locally well-posed for initial data in $u_{0} \in \dot{H}^{s_{c}}(\mathbf{R}^{3})$ and $u_{1} \in \dot{H}^{s_{c} - 1}(\mathbf{R}^{3})$ on some interval $[-T(u_{0}, u_{1}), T(u_{0}, u_{1})]$. The time of well-posedness $T(u_{0}, u_{1})$ depends on the profile of the initial data $(u_{0}, u_{1})$, not just its size. 

Additional regularity is enough to give a lower bound on the time of well-posedness. Therefore, there exists some $T(\| u_{0} \|_{\dot{H}^{s}}, \| u_{1} \|_{\dot{H}^{s - 1}}) > 0$ for any $s_{c} < s < \frac{3}{2}$.
\end{theorem}

There also is good reason to think that in the defocusing case the local solution ought to be global. In general, a solution to $(\ref{1.1})$ conserves the energy
\begin{equation}\label{1.5}
E(u(t)) = \frac{1}{2} \int u_{t}(t,x)^{2} dx + \frac{1}{2} \int |\nabla u(t,x)|^{2} dx + \frac{1}{p + 2} \int |u(t,x)|^{p + 1} dx.
\end{equation}
Since $(\ref{1.4})$ is energy-subcritical, conservation of energy implies that $(\ref{1.4})$ is globally well-posed for any initial data $u_{0} \in \dot{H}^{s_{c}} \cap \dot{H}^{1}$ and $u_{1} \in \dot{H}^{s_{c} - 1} \cap L^{2}$.\vspace{5mm}

Indeed, by the Sobolev embedding theorem,
\begin{equation}\label{1.6}
E(u(0)) \lesssim \| u_{t}(0) \|_{L^{2}(\mathbf{R}^{3})}^{2} + \| \nabla u(0) \|_{L^{2}(\mathbf{R}^{3})}^{2} + \| \nabla u(0) \|_{L^{2}(\mathbf{R}^{3})}^{2} \| u(0) \|_{\dot{H}^{s_{c}}(\mathbf{R}^{3})}^{p - 1},
\end{equation}
and therefore,
\begin{equation}\label{1.7}
E(u(0)) \lesssim_{\| u_{0} \|_{\dot{H}^{s_{c}}}} \| u_{t}(0) \|_{L^{2}}^{2} + \| \nabla u(0) \|_{L^{2}}^{2}.
\end{equation}
By conservation of energy, $E(u(0)) = E(u(t))$, and therefore $(\ref{1.7})$ gives a uniform bound over the norm $\| u_{t}(t) \|_{L^{2}}^{2} + \| \nabla u(t) \|_{L^{2}}^{2}$. Then since $(\ref{1.4})$ is energy-subcritical, a uniform bound over the energy is enough to ensure global well-posedness.\vspace{5mm}

\noindent \textbf{Remark:} This is not true for the focusing problem, which will not be discussed here.\vspace{5mm}

Moreover, the lack of a conserved quantity at the critical Sobolev regularity $s_{c} < 1$ is the only obstacle to proving global well-posedness and scattering for $(\ref{1.4})$ with radial data. Indeed,
\begin{theorem}\label{t1.2}
Suppose $u_{0} \in \dot{H}^{s_{c}}(\mathbf{R}^{3})$ and $u_{1} \in \dot{H}^{s_{c} - 1}(\mathbf{R}^{3})$ are radial functions, and $u$ solves $(\ref{1.4})$ on a maximal interval $0 \in I \subset \mathbf{R}$, with
\begin{equation}\label{1.8}
\sup_{t \in I} \| u(t) \|_{\dot{H}^{s_{c}}(\mathbf{R}^{3})} + \| u_{t}(t) \|_{\dot{H}^{s_{c} - 1}(\mathbf{R}^{3})} < \infty.
\end{equation}
Then $I = \mathbf{R}$ and the solution $u$ scatters both forward and backward in time.
\end{theorem}
\emph{Proof:} See \cite{Shen}. $\Box$\vspace{5mm}

The proof of Theorem $\ref{t1.2}$ in \cite{Shen} used the concentration compactness method. Such methods have been well utilized to study the quintic nonlinear wave equation
\begin{equation}\label{1.8.1}
u_{tt} - \Delta u + u^{5} = 0, \qquad u(0,x) = u_{0}, \qquad u_{t}(0,x) = u_{1}.
\end{equation}
The qualitative behavior of the quintic wave equation has been completely worked out, proving both global well-posedness and scattering, for both the radial (\cite{GSV}, \cite{Struwe}) and the nonradial case (\cite{BS}, \cite{Gril}, \cite{Shatah - Struwe}). The proof relies very heavily on conservation of the energy
\begin{equation}\label{1.8.2}
E(u(t)) = \frac{1}{2} \int u_{t}(t,x)^{2} dx + \frac{1}{2} \int |\nabla u(t,x)|^{2} dx + \frac{1}{6} \int u(t,x)^{6} dx.
\end{equation}
Conservation of energy ensures a uniform bound over the critical Sobolev norm, which guarantees that $(\ref{1.8})$ holds for $(\ref{1.8.1})$. Conservation of energy also yields a Morawetz estimate,
\begin{equation}\label{1.8.3}
\int \int \frac{u(t,x)^{6}}{|x|} dx dt \lesssim E(u(0)),
\end{equation}
which gives a space-time integral estimate for a solution to $(\ref{1.8.1})$.\vspace{5mm}

In order to make use of conservation of energy in the proof of Theorem $\ref{t1.1}$, the Fourier truncation method is used. The initial data is split into two pieces, a piece with small $\dot{H}^{s_{c}} \times \dot{H}^{s_{c} - 1}$ norm, and a piece with finite energy. Then, a solution $u$ to $(\ref{1.4})$ is shown to have the decomposition
\begin{equation}\label{1.9}
u(t) = v(t) + w(t),
\end{equation}
where $v(t)$ has uniformly bounded energy, and $w(t)$ is a small data scattering solution to $(\ref{1.4})$. By Theorem $\ref{t1.2}$, a uniform bound on the energy of $v(t)$ is enough to imply global well-posedness of $(\ref{1.4})$.\vspace{5mm}

\noindent \textbf{Remark:} The Fourier truncation method was used in \cite{KPV} to prove global well-posedness for the cubic problem when $s > \frac{3}{4}$.\vspace{5mm}

To prove scattering, the wave equation $(\ref{1.4})$ is rewritten in hyperbolic coordinates. These coordinates were quite useful to the cubic wave equation because the hyperbolic energy scales like the $\dot{H}^{1/2} \times \dot{H}^{-1/2}$ norm. For $3 < p < 5$, the hyperbolic energy and the energy ``sandwich" the $\dot{H}^{s_{c}} \times \dot{H}^{s_{c} - 1}$ norm, giving scattering.\vspace{5mm}

\noindent \textbf{Remark:} Previously, \cite{Shen1} used hyperbolic coordinates to prove scattering for $(\ref{1.4})$ with radial data lying in the energy space and a weighted Sobolev space. The weighted Sobolev space used in \cite{Shen1} also scales like the $\dot{H}^{1/2} \times \dot{H}^{-1/2}$ norm.\vspace{5mm}

As in \cite{D}, energy and hyperbolic energy bounds merely give a scattering size bound for any initial data in the critical Sobolev space, but with scattering size depending on the initial data $(u_{0}, u_{1})$, and not just its size. To prove a scattering size bound that depends on the size of the initial data, use Zorn's lemma. As in \cite{D2} and \cite{D}, it is shown by a profile decomposition that if $(u_{0}^{n}, u_{1}^{n}) \dot{H}^{s_{c}} \times \dot{H}^{s_{c} - 1}$ is a bounded sequence, then $\| u^{n} \|_{L_{t,x}^{p}(\mathbf{R} \times \mathbf{R}^{3})}$ is also uniformly bounded.\vspace{5mm}

\noindent \textbf{Remark:} The upper bound in $(\ref{1.4.1})$ is completely qualitative. Concentration compactness-type arguments that proved scattering in the energy-critical case also obtained a quantitative bound. See for example \cite{Tao}. Here we do not obtain any quantitative bounds at all. In the author's opinion, it would be very interesting to obtain some sort of quantitative bound.\vspace{5mm}

We begin by proving global well-posedness for the $p = 4$ case. We then generalize this global well-posedness result to any $3 < p < 5$. After proving global well-posedness, the hyperbolic coordinates are well-defined. In section four, we prove an estimate on the initial data, before obtaining a scattering bound in section five. We conclude with a Zorn's lemma argument in section six.\vspace{5mm}

\noindent \textbf{Acknowledgements:} The author was partially supported on NSF grant number $1764358$ during the writing of this paper. The author was also a guest of the Institute for Advanced Study during the writing of this paper.

\section{$p = 4$ case}
To simplify exposition by considering a specific case, consider $(\ref{1.1})$ with $p = 4$,
\begin{equation}\label{2.0}
u_{tt} - \Delta u + |u|^{3} u = 0.
\end{equation}
In this case
\begin{equation}\label{2.1}
s_{c} = \frac{3}{2} - \frac{2}{p} = \frac{5}{6}.
\end{equation}
Global well-posedness is proved by the Fourier truncation method. Split
\begin{equation}\label{2.1.1}
v_{0} = P_{\leq 1} u_{0}, \qquad w_{0} = P_{> 1} u_{0}, \qquad v_{1} = P_{\leq 1} u_{1}, \qquad w_{1} = P_{> 1} u_{1},
\end{equation}
and rescale by $(\ref{1.2})$ so that
\begin{equation}\label{2.2}
\| (w_{0}, w_{1}) \|_{\dot{H}^{5/6} \times \dot{H}^{-1/6}} < \epsilon.
\end{equation}
By Theorem $\ref{t1.2}$, $(\ref{2.0})$ has a local solution. Decompose the solution to $(\ref{2.0})$, $u = v + w$, where
\begin{equation}\label{2.3}
w_{tt} - \Delta w + |w|^{3} w = 0, \qquad w(0,x) = w_{0}, \qquad w_{t}(0,x) = w_{1},
\end{equation}
\begin{equation}\label{2.4}
v_{tt} - \Delta v + |u|^{3} u - |w|^{3} w = 0, \qquad v(0,x) = v_{0}, \qquad v_{t}(0,x) = v_{1}.
\end{equation}

Small data arguments and Strichartz estimates show that $(\ref{2.3})$ is globally well-posed and scattering.
\begin{theorem}\label{t2.1}
Let $I \subset \mathbf{R}$, $t_{0} \in I$, be an interval and let $u$ solve the linear wave equation
\begin{equation}\label{2.5}
u_{tt} - \Delta u = F, \hspace{5mm} u(t_{0}) = u_{0}, \hspace{5mm} u_{t}(t_{0}) = u_{1}.
\end{equation}
Then we have the estimates
\begin{equation}\label{2.6}
\aligned
\| u \|_{L_{t}^{p} L_{x}^{q}(I \times \mathbf{R}^{3})} + \| u \|_{L_{t}^{\infty} \dot{H}^{s}(I \times \mathbf{R}^{3})} + \| u_{t} \|_{L_{t}^{\infty} \dot{H}^{s - 1}(I \times \mathbf{R}^{3})} \\
\lesssim_{p, q, s, \tilde{p}, \tilde{q}} \| u_{0} \|_{\dot{H}^{s}(\mathbf{R}^{3})} + \| u_{1} \|_{\dot{H}^{s - 1}(\mathbf{R}^{3})} + \| F \|_{L_{t}^{\tilde{p}'} L_{x}^{\tilde{q}'}(I \times \mathbf{R}^{3})},
\endaligned
\end{equation}
whenever $s \geq 0$, $2 \leq p, \tilde{p} \leq \infty$, $2 \leq q, \tilde{q} < \infty$, and
\begin{equation}\label{2.7}
\frac{1}{p} + \frac{1}{q} \leq \frac{1}{2}, \hspace{5mm} \frac{1}{\tilde{p}} + \frac{1}{\tilde{q}} \leq \frac{1}{2}.
\end{equation}
\end{theorem}
\emph{Proof:} Theorem $\ref{t2.1}$ was proved for $p = q = 4$ in \cite{Stri} and then in \cite{GV} for a general choice of $(p, q)$. $\Box$\vspace{5mm}

Then,
\begin{equation}\label{2.8}
\| w \|_{L_{t,x}^{6} \cap L_{t}^{12/5} L_{x}^{12} \cap L_{t}^{\infty} \dot{H}^{5/6}} \lesssim \| (w_{0}, w_{1}) \|_{\dot{H}^{5/6} \times \dot{H}^{-1/6}} + \| w \|_{L_{t,x}^{6}}^{3} \| w \|_{L_{t}^{12/5} L_{x}^{12}} \lesssim \epsilon,
\end{equation}
which by $(\ref{2.2})$ implies that $w$ is scattering.

Also, by the radial Strichartz estimate and Bernstein's inequality,
\begin{equation}\label{2.4.2}
\| w \|_{L_{t}^{\infty} L_{x}^{2}} \lesssim \epsilon + \| w \|_{L_{t,x}^{6}}^{3} \| w \|_{L_{t}^{\infty} L_{x}^{2}} \lesssim \epsilon.
\end{equation}
\begin{theorem}[Radial Strichartz estimate]\label{t2.2}
For $(u_{0}, u_{1})$ radially symmetric, and $u$ solves $(\ref{2.5})$ with $F = 0$,
\begin{equation}\label{2.9}
\| u \|_{L_{t}^{2} L_{x}^{\infty}(\mathbf{R} \times \mathbf{R}^{3})} \lesssim \| u_{0} \|_{\dot{H}^{1}(\mathbf{R}^{3})} + \| u_{1} \|_{L^{2}(\mathbf{R}^{3})}.
\end{equation}
\end{theorem}
\emph{Proof:} This theorem was proved in \cite{KlMa}. The dual of $(\ref{2.9})$ is that if $u_{0} = u_{1} = 0$, and $F$ is radial, then
\begin{equation}\label{2.9.1}
\| u \|_{L_{t}^{\infty} L_{x}^{2}} \lesssim \| F \|_{L_{t}^{2} L_{x}^{1}}.
\end{equation}
$\Box$\vspace{5mm}

Next, let $E(t)$ be the energy of $v$,
\begin{equation}\label{2.10}
E(t) = \frac{1}{2} \int |\nabla v|^{2} + \frac{1}{2} \int v_{t}^{2} + \frac{1}{5} \int |v|^{5} dx.
\end{equation}
By the Sobolev embedding theorem and $(\ref{2.2})$,
\begin{equation}\label{2.11}
E(0) \lesssim (\| u_{0} \|_{\dot{H}^{5/6}} + \| u_{1} \|_{\dot{H}^{-1/6}})^{2} + (\| u_{0} \|_{\dot{H}^{5/6}} + \| u_{1} \|_{\dot{H}^{-1/6}})^{5}.
\end{equation}
To prove global well--posedness it is enough to prove a uniform bound on $E(t)$.
\begin{theorem}\label{t2.3}
The energy $E(t)$ given by $(\ref{2.10})$ is uniformly bounded for all $t \in \mathbf{R}$, and moreover,
\begin{equation}\label{2.12}
\sup_{t \in \mathbf{R}} E(t) \lesssim_{\| u_{0} \|_{\dot{H}^{5/6}}, \| u_{1} \|_{\dot{H}^{-1/6}}} E(0).
\end{equation}
\end{theorem}
\emph{Proof:} The proof is quite similar to the proof in \cite{D}. By direct computation,
\begin{equation}\label{2.13}
\frac{d}{dt} E(v(t)) = \int v_{t} [|v + w|^{3} (v + w) - |w|^{3} w - |v|^{3} v] dx.
\end{equation}
By Taylor's theorem,
\begin{equation}\label{2.14}
\aligned
|v + w|^{3} (v + w) - |v|^{3} v - |w|^{3} w = 4w \int_{0}^{1} |v + \tau w|^{3} d\tau - 4w \int_{0}^{1} |\tau w|^{3} d\tau \\
= 12 wv \int_{0}^{1} \int_{0}^{1} |sv + \tau w|(sv + \tau w) ds d\tau = 4 |v|^{3} w + O(|v|^{2} |w|^{2}) + O(|v| |w|^{3}).
\endaligned
\end{equation}
By H{\"o}lder's inequality and $(\ref{2.10})$,
\begin{equation}\label{2.15}
\langle v_{t}, |v|^{2} |w|^{2} \rangle \lesssim \| v_{t} \|_{L_{x}^{2}(\mathbf{R}^{3})} \| v \|_{L_{x}^{6}(\mathbf{R}^{3})}^{1/3} \| v \|_{L_{x}^{5}(\mathbf{R}^{3})}^{5/3} \| w \|_{L_{x}^{18}(\mathbf{R}^{3})}^{2} \lesssim E(t) \| w(t) \|_{L_{x}^{18}(\mathbf{R}^{3})}^{2},
\end{equation}
and
\begin{equation}\label{2.16}
\langle v_{t}, |v| |w|^{3} \rangle \lesssim \| v_{t} \|_{L_{x}^{2}(\mathbf{R}^{3})} \| v \|_{L_{x}^{6}(\mathbf{R}^{3})} \| w \|_{L_{x}^{9}(\mathbf{R}^{3})}^{3} \lesssim E(t) \| w(t) \|_{L_{x}^{9}(\mathbf{R}^{3})}^{3}.
\end{equation}
Therefore,
\begin{equation}\label{2.17}
\frac{d}{dt} E(t) = 4 \langle v_{t}, |v|^{3} w \rangle + E(t) [\| w(t) \|_{L_{x}^{18}(\mathbf{R}^{3})}^{2} + \| w(t) \|_{L_{x}^{9}(\mathbf{R}^{3})}^{3}].
\end{equation}

If the term $4 \langle v_{t}, |v|^{3} w \rangle$ could be dropped, and
\begin{equation}\label{2.18}
\frac{d}{dt} E(t) \lesssim E(t) [\| w(t) \|_{L_{x}^{18}(\mathbf{R}^{3})}^{2} + \| w(t) \|_{L_{x}^{9}(\mathbf{R}^{3})}^{3}],
\end{equation}
then by radial Strichartz estimates, $(\ref{2.1.1})$, and $(\ref{2.8})$,
\begin{equation}\label{2.19}
\int_{\mathbf{R}} \| w(t) \|_{L_{x}^{18}(\mathbf{R}^{3})}^{2} + \| w(t) \|_{L_{x}^{9}(\mathbf{R}^{3})}^{3} dt \lesssim \epsilon^{2},
\end{equation}
which by Gronwall's inequality implies $\sup_{t \in \mathbf{R}} E(t) \lesssim E(0)$.
\begin{theorem}[Radial Strichartz estimates]\label{t2.3}
Let $(u_{0}, u_{1})$ be spherically symmetric, and suppose $u$ solves $(\ref{2.5})$ with $F = 0$. Then if $q > 4$ and
\begin{equation}\label{2.20}
\frac{1}{2} + \frac{3}{q} = \frac{3}{2} - s,
\end{equation}
then
\begin{equation}\label{2.21}
\| u \|_{L_{t}^{2} L_{x}^{q}(\mathbf{R} \times \mathbf{R}^{3})} \lesssim \| u_{0} \|_{\dot{H}^{s}(\mathbf{R}^{3})} + \| u_{1} \|_{\dot{H}^{s - 1}(\mathbf{R}^{3})}.
\end{equation}
\end{theorem}
\emph{Proof:} This was proved in \cite{Sterb}. $\Box$\vspace{5mm}

As in \cite{D}, the contribution of $4 \langle v_{t}, |v|^{3} w \rangle$ will be controlled by replacing $E(t)$ with a term $\mathcal E(t) \sim E(t)$ that has better time differentiability properties. Define
\begin{equation}\label{2.22}
\mathcal E(t) = E(t) + c M(t) - \int |v|^{3} v w dx,
\end{equation}
where $M(t)$ is the Morawetz potential
\begin{equation}\label{2.23}
M(t) = \int v_{t} \frac{x}{|x|} \cdot \nabla v dx + \int v_{t} \frac{1}{|x|} v dx,
\end{equation}
and $c > 0$ is a small, fixed constant. By Hardy's inequality,
\begin{equation}\label{2.24}
M(t) \lesssim c \| \nabla v \|_{L^{2}(\mathbf{R}^{3})} \| v_{t} \|_{L^{2}(\mathbf{R}^{3})} \lesssim c E(t),
\end{equation}
and by $(\ref{2.8})$,
\begin{equation}\label{2.25}
\aligned
\int |v|^{3} v w dx  \lesssim \| v \|_{L_{x}^{5}}^{10/3} \| v \|_{L_{x}^{6}}^{2/3} \| w \|_{L_{x}^{9/2}} \lesssim \epsilon E(t).
\endaligned
\end{equation}
Therefore, $\mathcal E(t) \sim E(t)$.

Next, by the product rule,
\begin{equation}\label{2.26}
4 \langle v_{t}, |v|^{3} w \rangle - \frac{d}{dt} \int |v|^{3} v w dx = \langle v, |v|^{3} w_{t} \rangle.
\end{equation}
Also, by direct computation and integrating by parts, since $v$ is radial,
\begin{equation}\label{2.27}
\aligned
c \frac{d}{dt} M(t) = -\frac{c}{2} v(t,0)^{2} - \frac{3c}{5} \int \frac{|v(t,x)|^{5}}{|x|} dx \\ - c \int (|v + w|^{3}(v + w) - |v|^{3} v - |w|^{3} w) \frac{x}{|x|} \cdot \nabla v dx \\
- c \int (|v + w|^{3}(v + w) - |v|^{3} v - |w|^{3} w) \frac{1}{|x|} v dx.
\endaligned
\end{equation}
\textbf{Remark:} The virial identities will be computed in more detail in the next section.

Therefore,
\begin{equation}\label{2.28}
\aligned
\frac{d}{dt} \mathcal E(t) = -\frac{c}{2} v(t,0)^{2} - \frac{3c}{5} \int \frac{|v(t,x)|^{5}}{|x|} dx \\ - c \int (|v + w|^{3}(v + w) - |v|^{3} v - |w|^{3} w) \frac{x}{|x|} \cdot \nabla v dx + \langle v, |v|^{3} w_{t} \rangle \\
- c \int (|v + w|^{3}(v + w) - |v|^{3} v - |w|^{3} w) \frac{1}{|x|} v dx \\ 
+ O(E(t) [\| w(t) \|_{L_{x}^{18}(\mathbf{R}^{3})}^{2} + \| w(t) \|_{L_{x}^{9}(\mathbf{R}^{3})}^{3}]).
\endaligned
\end{equation}

By Hardy's inequality, the Sobolev embedding theorem, and the Cauchy--Schwartz inequality,
\begin{equation}\label{2.29}
\aligned
\int (|v + w|^{3}(v + w) - |v|^{3} v - |w|^{3} w) \frac{1}{|x|} v dx \\
\lesssim (\int \frac{1}{|x|} v^{5} dx)^{2/3} \cdot \| \frac{1}{|x|^{1/2}} v \|_{L_{x}^{3}}^{2/3} \| w \|_{L_{x}^{9}} + \| \frac{1}{|x|} v \|_{L^{2}} \| v \|_{L_{x}^{6}} \| w \|_{L_{x}^{9}}^{3} \\
\lesssim \delta (\int \frac{1}{|x|} |v|^{5} dx) + \frac{1}{\delta} E(t) \| w(t) \|_{L_{x}^{9}}^{3}.
\endaligned
\end{equation}
Also, following $(\ref{2.15})$ and $(\ref{2.16})$,
\begin{equation}\label{2.30}
\aligned
c \int [|v|^{2} |w|^{2} + |v| |w|^{3}] \frac{x}{|x|} \cdot \nabla v dx \lesssim \| \nabla v \|_{L^{2}} \| v \|_{L^{6}}^{1/3} \| v \|_{L^{5}}^{5/3} \| w \|_{L_{x}^{18}}^{2} \\ + \| \nabla v \|_{L_{x}^{2}} \| v \|_{L_{x}^{6}} \| w \|_{L_{x}^{9}}^{3} \lesssim E(t) [\| w \|_{L_{x}^{18}}^{2} + \| w \|_{L_{x}^{9}}^{3}].
\endaligned
\end{equation}
Therefore,
\begin{equation}\label{2.31}
\aligned
\frac{d}{dt} \mathcal E(t) + \frac{c}{2} v(t,0)^{2} + \frac{3c}{5} \int \frac{|v(t,x)|^{5}}{|x|} dx \\ + c \int w \frac{x}{|x|} \cdot \nabla (|v|^{3} v) dx + \langle v, |v|^{3} w_{t} \rangle \\
\lesssim \frac{1}{\delta} E(t) [\| w(t) \|_{L_{x}^{18}(\mathbf{R}^{3})}^{2} + \| w(t) \|_{L_{x}^{9}(\mathbf{R}^{3})}^{3}] + \delta (\int \frac{|v(t,x)|^{5}}{|x|} dx).
\endaligned
\end{equation}
Also, by Lemma $\ref{l3.2}$, if $P_{j}$ is a Littlewood--Paley projection operator,
\begin{equation}\label{2.32}
\int \frac{1}{|x|} |P_{\leq j} v|^{5} dx + \int \frac{1}{|x|} |P_{\geq j} v|^{5} dx \lesssim \int \frac{1}{|x|} |v|^{5} dx.
\end{equation}

Making a Littlewood--Paley decomposition,
\begin{equation}\label{2.33}
\langle |v|^{3} v, w_{t} \rangle = \sum_{j} \langle |v|^{3} v, P_{j} w_{t} \rangle.
\end{equation}
Then by H{\"o}lder's inequality, $(\ref{2.32})$, and the Cauchy--Schwartz inequality,
\begin{equation}\label{2.34}
\aligned
\sum_{j} \langle |v|^{3} v - |P_{\leq j} v|^{3} (P_{\leq j} v), P_{j} w_{t} \rangle \\ \lesssim \sum_{j} \| |x|^{1/10} P_{\geq j} v \|_{L_{x}^{5/2}} (\int \frac{1}{|x|} (|P_{\leq j} v|^{5} + |P_{\geq j} v|^{5}) dx)^{3/5} \| |x|^{1/2} P_{j} w_{t} \|_{L_{x}^{\infty}} \\
\lesssim (\int \frac{1}{|x|} |v|^{5} dx)^{3/5} \cdot \sum_{j} \| |x|^{1/10} P_{\geq j} v \|_{L_{x}^{5/2}} \| |x|^{1/2} P_{j} w_{t} \|_{L_{x}^{\infty}}.
\endaligned
\end{equation}
By Bernstein's inequality and the radial Sobolev embedding theorem,
\begin{equation}\label{2.35}
\| |x|^{1/10} P_{\geq j} v \|_{L_{x}^{5/2}(\mathbf{R}^{3})} \lesssim 2^{-4j/5} \| \nabla v \|_{L_{x}^{2}(\mathbf{R}^{3})} \lesssim 2^{-4j/5} E(t)^{1/2}.
\end{equation}
Also, by Bernstein's inequality and the radial Sobolev embedding theorem,
\begin{equation}\label{2.36}
\aligned
\langle |P_{\leq j} v|^{3} (P_{\leq j} v), P_{j} w_{t} \rangle \\ \lesssim 2^{-j} \| |x|^{1/10} \nabla P_{\leq j} v \|_{L_{x}^{5/2}(\mathbf{R}^{3})} (\int \frac{1}{|x|} |P_{\leq j} v|^{5} dx)^{3/5} \| |x|^{1/2} P_{j} w_{t} \|_{L_{x}^{\infty}(\mathbf{R}^{3})} \\
\lesssim 2^{-4j/5} E(t)^{1/2} (\int \frac{1}{|x|} |v|^{5} dx)^{3/5} \| |x|^{1/2} P_{j} w_{t} \|_{L_{x}^{\infty}}.
\endaligned
\end{equation}
Similarly, by Bernstein's inequality and $(\ref{2.32})$,
\begin{equation}\label{2.37}
\aligned
\int (P_{j} w) \frac{x}{|x|} \cdot \nabla(|P_{\leq j} v|^{3} (P_{\leq j} v)) dx \\ \lesssim 2^{-j} \| |x|^{1/2} P_{j} \nabla w \|_{L_{x}^{\infty}} \| |x|^{1/10} \nabla P_{\leq j} v \|_{L_{x}^{5/2}} (\frac{1}{|x|} |P_{\leq j} v|^{5} dx)^{3/5} \\
\lesssim 2^{-4j/5} E(t)^{1/2} (\int \frac{1}{|x|} |v|^{5} dx)^{3/5} \| |x|^{1/2} P_{j} \nabla w \|_{L_{x}^{\infty}}.
\endaligned
\end{equation}
Meanwhile, integrating by parts,
\begin{equation}\label{2.38}
\aligned
\int (P_{j} w) \frac{x}{|x|} \cdot \nabla(|v|^{3} v - |P_{\leq j} v|^{3} (P_{\leq j} v)) dx \\
= -\int (P_{j} \nabla w) \cdot \frac{x}{|x|} (|v|^{3} v - |P_{\leq j} v|^{3} (P_{\leq j} v)) dx \\
-2 \int (P_{j} w) \frac{1}{|x|} (|v|^{3} v - |P_{\leq j} v|^{3} (P_{\leq j} v)) dx.
\endaligned
\end{equation}
The term
\begin{equation}\label{2.39}
-\int (P_{j} \nabla w) \cdot \frac{x}{|x|} (|v|^{3} v - |P_{\leq j} v|^{3} (P_{\leq j} v)) dx
\end{equation}
may be handled in a manner identical to $(\ref{2.36})$, giving
\begin{equation}\label{2.40}
(\ref{2.39}) \lesssim 2^{-4j/5} E(t)^{1/2} (\int \frac{1}{|x|} |v|^{5} dx)^{3/5} \| |x|^{1/2} P_{j} \nabla w \|_{L_{x}^{\infty}}.
\end{equation}
Meanwhile, by $(\ref{2.29})$ and $(\ref{2.32})$,
\begin{equation}\label{2.41}
-2 \int (P_{j} w) \frac{1}{|x|} (|v|^{3} v - |P_{\leq j} v|^{3} (P_{\leq j} v)) dx \lesssim \delta (\int \frac{1}{|x|} |v|^{5} dx) + \frac{1}{\delta} E(t) \| w(t) \|_{L_{x}^{9}}^{3}.
\end{equation}

Therefore, by $(\ref{2.31})$--$(\ref{2.41})$,
\begin{equation}\label{2.42}
\aligned
\frac{d}{dt} \mathcal E(t) + \frac{c}{2} v(t,0)^{2} + \frac{3c}{5} \int \frac{|v(t,x)|^{5}}{|x|} dx \\ 
\lesssim \frac{1}{\delta} E(t) [\| w(t) \|_{L_{x}^{18}(\mathbf{R}^{3})}^{2} + \| w(t) \|_{L_{x}^{9}(\mathbf{R}^{3})}^{3}] + \delta (\int \frac{|v(t,x)|^{5}}{|x|} dx) \\ + E(t)^{5/4} (\sum_{j} 2^{-4j/5} \| |x|^{1/2} P_{j} \nabla_{t,x} w \|_{L_{x}^{\infty}})^{5/2}.
\endaligned
\end{equation}
For $\delta > 0$ sufficiently small,
\begin{equation}\label{2.43}
\delta (\int \frac{1}{|x|} |v|^{5} dx)
\end{equation}
may be absorbed into the left hand side of $(\ref{2.42})$.

Using Corollary $3.3$ from \cite{D}, by Bernstein's inequality and $(\ref{2.8})$
\begin{equation}\label{2.44}
\sum_{j \geq 0} 2^{-4j/5} \| |x|^{1/2} P_{j} \nabla_{t,x} w \|_{L_{t}^{5/2} L_{x}^{\infty}(\mathbf{R} \times \mathbf{R}^{3})} \lesssim \| w_{0} \|_{\dot{H}^{5/6}(\mathbf{R}^{3})} + \| w_{1} \|_{\dot{H}^{-1/6}(\mathbf{R}^{3})} \lesssim \epsilon.
\end{equation}
Also by Bernstein's inequality and $(\ref{2.4.2})$,
\begin{equation}\label{2.45}
\sum_{j \leq 0} 2^{-4j/5} \| |x|^{1/2} P_{j} \nabla_{t,x} w \|_{L_{t}^{5/2} L_{x}^{\infty}(\mathbf{R} \times \mathbf{R}^{3})} \lesssim \| w_{0} \|_{L^{2}(\mathbf{R}^{3})} + \| w_{1} \|_{\dot{H}^{-1}(\mathbf{R}^{3})} \lesssim \epsilon.
\end{equation}
Therefore, by $(\ref{2.11})$, $(\ref{2.19})$, $(\ref{2.44})$, $(\ref{2.45})$, and Gronwall's inequality, for $\epsilon(\| u_{0} \|_{\dot{H}^{5/6}}, \| u_{1} \|_{\dot{H}^{-1/6}})$ sufficiently small, $(\ref{2.12})$ holds, proving Theorem $\ref{t2.3}$. $\Box$

\section{Global well-posedness for general $p$}
Proof of global well-posedness for a general $(\ref{1.1})$ is a generalization of the $p = 4$ case.
\begin{theorem}\label{t3.0}
The nonlinear wave equation
\begin{equation}\label{3.1}
u_{tt} - \Delta + |u|^{p - 1} u = 0, \qquad u(0,x) = u_{0}, \qquad u_{t}(0,x) = u_{1},
\end{equation}
with radial initial data $u_{0} \in \dot{H}^{s_{c}}(\mathbf{R}^{3})$, $u_{1} \in \dot{H}^{s_{c} - 1}(\mathbf{R}^{3})$, with $s_{c} = \frac{3}{2} - \frac{2}{p - 1}$, $3 < p < 5$, is globally well-posed.
\end{theorem}
\emph{Proof:} The proof is a generalization of the argument in the $p = 4$ case.\vspace{5mm}

First prove a generalized Morawetz inequality.
\begin{theorem}[Morawetz inequality]\label{t3.1}
If $u$ solves $(\ref{3.1})$ on an interval $I$, then
\begin{equation}\label{3.2}
\int_{I} \int \frac{|u(t,x)|^{p + 1}}{|x|} dx dt \lesssim E(u),
\end{equation}
where $E$ is the conserved energy $(\ref{1.5})$.
\end{theorem}
\emph{Proof:} Define the Morawetz potential
\begin{equation}\label{3.3}
M(t) = \int u_{t} u_{r} r^{2} dr + \int u_{t} u r dr.
\end{equation}
By direct computation,
\begin{equation}\label{3.4}
\frac{d}{dt} M(t) = -\frac{1}{2} u(t,0)^{2} - \frac{p - 1}{p + 1} \int |u|^{p + 1} r dr.
\end{equation}
Then $(\ref{3.2})$ holds by the fundamental theorem of calculus and Hardy's inequality. $\Box$\vspace{5mm}

The Morawetz estimate commutes very well with Littlewood--Paley projections.
\begin{lemma}\label{l3.2}
For any $j$,
\begin{equation}\label{3.5}
\int \frac{1}{|x|} |P_{\leq j} v|^{p + 1} dx + \int \frac{1}{|x|} |P_{\geq j} v|^{p + 1} dx \lesssim \int \frac{1}{|x|} |v|^{p + 1} dx.
\end{equation}
\end{lemma}
\noindent \emph{Proof:} Let $\psi$ be the Littlewood--Paley kernel.
\begin{equation}\label{3.6}
\frac{1}{|x|^{\frac{1}{p + 1}}} P_{\leq j} v(x) = \frac{1}{|x|^{\frac{1}{p + 1}}} \int 2^{3j} \psi(2^{j}(x - y)) v(y) dy.
\end{equation}
When $|y| \lesssim |x|$,
\begin{equation}\label{3.7}
\frac{1}{|x|^{\frac{1}{p + 1}}} 2^{3j} \psi(2^{j}(x - y)) \lesssim 2^{3j} \psi(2^{j}(x - y)) \frac{1}{|y|^{\frac{1}{p + 1}}}.
\end{equation}
When $|y| \gg |x|$ and $|x| \geq 2^{-j}$, since $\psi$ is rapidly decreasing, for any $N$,
\begin{equation}\label{3.8}
\aligned
\frac{1}{|x|^{\frac{1}{p + 1}}} 2^{3j} \psi(2^{j}(x - y)) \lesssim_{N} \frac{1}{|x|^{\frac{1}{p + 1}}} \frac{2^{3j}}{(1 + 2^{j} |x - y|)^{N}} \\ \lesssim \frac{1}{|x|^{\frac{1}{p + 1}} 2^{j} |y|} \frac{2^{3j}}{(1 + 2^{j} |x - y|)^{N - 1}} \lesssim \frac{1}{|y|^{\frac{1}{p + 1}}} \frac{2^{3j}}{(1 + 2^{j}|x - y|)^{N - 1}}.
\endaligned
\end{equation}
Combining $(\ref{3.7})$ and $(\ref{3.8})$,
\begin{equation}\label{3.9}
\| \frac{1}{|x|^{\frac{1}{p + 1}}} |P_{\leq j} v| \|_{L^{p + 1}(|x| \geq 2^{-j})} \lesssim \| \frac{1}{|x|^{\frac{1}{p + 1}}} v \|_{L^{p + 1}(\mathbf{R}^{3})}.
\end{equation}
When $|y| \gg |x|$ and $|x| \leq 2^{-j}$, since $\psi$ is rapidly decreasing, for any $N$,
\begin{equation}\label{3.10}
\aligned
\frac{1}{|x|^{\frac{1}{p + 1}}} 2^{3j} \psi(2^{j}(x - y)) \lesssim_{N} \frac{1}{|x|^{\frac{1}{p + 1}}} \frac{2^{3j}}{(1 + 2^{j} |x - y|)^{N}} \\ \lesssim \frac{1}{|x|^{\frac{1}{p + 1}}} \frac{2^{3j}}{(1 + 2^{j} |x - y|)^{N - \frac{1}{p + 1}}} \frac{1}{2^{\frac{j}{p + 1}} |y|^{\frac{1}{p + 1}}}.
\endaligned
\end{equation}
\begin{equation}\label{3.11}
\| \frac{2^{3j - \frac{j}{p + 1}}}{(1 + 2^{j}|x - y|)^{N}} \|_{L^{\frac{p + 1}{p}}(\mathbf{R}^{3})} \lesssim 2^{\frac{2j}{p + 1}},
\end{equation}
so by $(\ref{3.8})$, Young's inequality, and H{\"o}lder's inequality,
\begin{equation}\label{3.12}
\| \frac{1}{|x|^{\frac{1}{p + 1}}} |P_{\leq j} v| \|_{L^{p + 1}(|x| \leq 2^{-j})} \lesssim \| \frac{1}{|x|^{\frac{1}{p + 1}}} v \|_{L^{p + 1}(\mathbf{R}^{3})}.
\end{equation}
This proves $(\ref{3.5})$. $\Box$\vspace{5mm}

Next, split a local solution $(\ref{3.1})$, $u = v + w$, where $w$ solves
\begin{equation}\label{3.13}
w_{tt} - \Delta w + |w|^{p - 1} w = 0, \qquad w(0,x) = w_{0}, \qquad w_{t}(0,x) = w_{1},
\end{equation}
and $v$ solves
\begin{equation}\label{3.14}
v_{tt} - \Delta v + |u|^{p - 1} u - |w|^{p - 1} w = 0, \qquad v(0,x) = v_{0}, \qquad v_{t}(0,x) = v_{1}.
\end{equation}
As usual, use the rescaling $(\ref{1.2})$ so that $v_{0} = P_{\leq 1} u_{0}$, $v_{1} = P_{\leq 1} u_{1}$, $w_{0} = P_{> 1} u_{0}$, $w_{1} = P_{> 1} u_{1}$, and
\begin{equation}\label{3.15}
\| w_{0} \|_{\dot{H}^{s_{c}}(\mathbf{R}^{3})} + \| w_{1} \|_{\dot{H}^{s_{c} - 1}(\mathbf{R}^{3})} < \epsilon.
\end{equation}
As in $(\ref{2.11})$,
\begin{equation}\label{3.16}
E(0) \lesssim (\| u_{0} \|_{\dot{H}^{s_{c}}} + \| u_{1} \|_{\dot{H}^{s_{c} - 1}})^{2} + (\| u_{0} \|_{\dot{H}^{s_{c}}} + \| u_{1} \|_{\dot{H}^{s_{c} - 1}})^{p + 1}.
\end{equation}

By small data arguments, $(\ref{3.13})$ is globally well-posed and scattering for $\epsilon > 0$ sufficiently small. Indeed,
\begin{equation}\label{3.17}
\| w \|_{L_{t,x}^{2(p - 1)} \cap L_{t}^{\frac{2}{s_{c}}} L_{x}^{\frac{2}{1 - s_{c}}}} \lesssim \| w_{0} \|_{\dot{H}^{s_{c}}} + \| w_{1} \|_{\dot{H}^{s_{c} - 1}} + \| w \|_{L_{t,x}^{2(p - 1)}}^{p - 1} \| w \|_{L_{t}^{\frac{2}{s_{c}}} L_{x}^{\frac{2}{1 - s_{c}}}} < \epsilon,
\end{equation}
and by Bernstein's inequality,
\begin{equation}\label{3.17.1}
\| w \|_{L_{t}^{\infty} L_{x}^{2}} \lesssim \| w_{0} \|_{L^{2}} + \| w_{1} \|_{\dot{H}^{-1}} + \| w \|_{L_{t}^{\infty} L_{x}^{2}} \| w \|_{L_{t,x}^{2(p - 1)}}^{p - 1} < \epsilon.
\end{equation}

Now define the energy of $v$,
\begin{equation}\label{3.18}
E(t) = \frac{1}{2} \int |\nabla v|^{2} + \frac{1}{2} \int v_{t}(t,x)^{2} dx + \frac{1}{p + 1} \int |v(t,x)|^{p + 1} dx,
\end{equation}
and let
\begin{equation}\label{3.19}
\mathcal E(t) = E(t) + c M(t) - \int |v|^{p - 1} vw dx,
\end{equation}
where $c > 0$ is a small constant and $M(t)$ is given by $(\ref{3.3})$, with $u$ replaced by $v$. Then by $(\ref{2.13})$ and $(\ref{3.4})$,
\begin{equation}\label{3.20}
\aligned
\frac{d}{dt} \mathcal E(t) + \frac{c}{2} v(t,0)^{2} + c(1 - \frac{2}{p + 1}) \int \frac{|v(t,x)|^{p + 1}}{|x|} dx \\ = - \langle v_{t}, |v + w|^{p - 1} (v + w) - |v|^{p - 1} v - |w|^{p - 1} w \rangle + \frac{d}{dt} \int |v|^{p - 1} vw dx \\
- c \int [|v + w|^{p - 1} (v + w) - |v|^{p - 1} v - |w|^{p - 1} w] \frac{x}{|x|} \cdot \nabla v dx \\
- c \int [|v + w|^{p - 1} (v + w) - |v|^{p - 1} v - |w|^{p - 1} w] \frac{1}{|x|} v dx.
\endaligned
\end{equation}
By $(\ref{2.29})$, Hardy's inequality, and the Cauchy--Schwartz inequality,
\begin{equation}\label{3.21}
\aligned
\int [|v + w|^{p - 1} (v + w) - |v|^{p - 1} v - |w|^{p - 1} w] \frac{1}{|x|} v dx \\ \lesssim (\int \frac{1}{|x|} |v|^{p + 1} dx)^{\frac{p - 2}{p - 1}} \| \frac{1}{|x|^{1/2}} v \|_{L_{x}^{3}}^{\frac{2}{p - 1}} \| w \|_{L^{3(p - 1)}} + \| \frac{1}{|x|} v \|_{L^{2}} \| v \|_{L^{6}} \| w \|_{L^{3(p - 1)}}^{p - 1} \\ \lesssim \delta (\int \frac{1}{|x|} |v|^{p + 1} dx) + \frac{1}{\delta} E(t) \| w \|_{L^{3(p - 1)}}^{p - 1}.
\endaligned
\end{equation}
Also by $(\ref{2.14})$,
\begin{equation}\label{3.22}
|v + w|^{p - 1} (v + w) - |v|^{p - 1} v - |w|^{p - 1} w = p |v|^{p - 1} w + O(|v|^{p - 2} |w|^{2}) + O(|v| |w|^{p - 1}).
\end{equation}
Then by H{\"o}lder's inequality,
\begin{equation}\label{3.23}
\aligned
\int [O(|v|^{p - 2} |w|^{2}) + O(|v| |w|^{p - 1})] \frac{x}{|x|} \cdot \nabla v dx \\ \lesssim \| \nabla v \|_{L^{2}} \| v \|_{L^{6}} \| w \|_{L^{3(p - 1)}}^{p - 1} + E(t) \| w \|_{L^{\frac{3}{1 - s_{c}}}}^{2} \\
\lesssim E(t) [\| w \|_{L^{3(p - 1)}}^{p - 1} + \| w \|_{L^{\frac{3}{1 - s_{c}}}}^{2}],
\endaligned
\end{equation}
and
\begin{equation}\label{3.24}
\langle v_{t}, [O(|v|^{p - 2} |w|^{2}) + O(|v| |w|^{p - 1})] \rangle \lesssim E(t) [\| w \|_{L^{3(p - 1)}}^{p - 1} + \| w \|_{L^{\frac{3}{1 - s_{c}}}}^{2}].
\end{equation}
Now then, by the product rule,
\begin{equation}\label{3.25}
p \langle v_{t}, |v|^{p - 1} w \rangle - \frac{d}{dt} \int |v|^{p - 1} v w dx = \langle |v|^{p - 1} v, w_{t} \rangle.
\end{equation}
Following $(\ref{2.34})$,
\begin{equation}\label{3.26}
\aligned
\sum_{j} \langle |v|^{p - 1} v - |P_{\leq j} v|^{p - 1} (P_{\leq j} v), P_{j} w_{t} \rangle \\ \lesssim (\int \frac{1}{|x|} |P_{\leq j} v|^{p + 1} + \frac{1}{|x|} |P_{\geq j} v|^{p + 1})^{\frac{p - 1}{p + 1}} \sum_{j} \| |x|^{1/2} P_{j} w_{t} \|_{L^{\infty}} \| |x|^{\frac{p - 3}{2(p + 1)}} |P_{\geq j} v| \|_{L^{\frac{p + 1}{2}}} \\
\lesssim \sum_{j} 2^{-\frac{4j}{p + 1}} E(t)^{1/2} (\int \frac{1}{|x|} |v|^{p + 1} dx)^{\frac{p - 1}{p + 1}} \| P_{j} w_{t} \|_{L^{\infty}} \\
\lesssim \delta (\int \frac{1}{|x|} |v|^{p + 1} dx) + \frac{1}{\delta} E(t)^{\frac{p + 1}{4}} (\sum_{j} 2^{-\frac{4j}{p + 1}} \| P_{j} w_{t} \|_{L^{\infty}})^{\frac{p + 1}{4}}.
\endaligned
\end{equation}
Integrating by parts as in $(\ref{2.37})$,
\begin{equation}
\sum_{j} \langle |P_{\leq j} v|^{p - 1} (P_{\leq j} v), P_{j} w_{t} \rangle \lesssim \delta (\int \frac{1}{|x|} |v|^{p + 1} dx) + \frac{1}{\delta} E(t)^{\frac{p + 1}{4}} (\sum_{j} 2^{-\frac{4j}{p + 1}} \| P_{j} w_{t} \|_{L^{\infty}})^{\frac{p + 1}{4}}.
\end{equation}
By $(\ref{3.17})$ and $(\ref{3.17.1})$,
\begin{equation}\label{3.27}
E(t)^{\frac{p + 1}{4}} \int_{\mathbf{R}} (\sum_{j} 2^{-\frac{4j}{p + 1}} \| P_{j} w_{t} \|_{L^{\infty}})^{\frac{p + 1}{4}} dt \lesssim E(t)^{\frac{p + 1}{4}} E(0)^{\frac{3 - p}{4}} \epsilon.
\end{equation}
Then by Gronwall's inequality, $(\ref{3.24})$, and $(\ref{3.27})$,
\begin{equation}\label{3.28}
\sup_{t \in \mathbf{R}} \mathcal E(t) \lesssim \mathcal E(0),
\end{equation}
which completes the proof of Theorem $\ref{t3.0}$. $\Box$

\section{Scattering: Estimates on initial data}
To prove scattering, let $\phi(x)$ be a smooth function supported on $|x| \leq 1$ and $\phi(x) = 1$ on $|x| \leq \frac{1}{2}$. Then for $R > 0$ sufficiently large,
\begin{equation}\label{4.1}
\| (1 - \phi(\frac{x}{R})) u_{0} \|_{\dot{H}^{s_{c}}(\mathbf{R}^{3})} + \| (1 - \phi(\frac{x}{R})) u_{1} \|_{\dot{H}^{s_{c} - 1}(\mathbf{R}^{3})} < \epsilon.
\end{equation}
Then rescale according to $(\ref{1.2})$,
\begin{equation}\label{4.2}
u_{0}(x) \mapsto (2R)^{\frac{2}{p - 1}} u_{0}(2Rx), \qquad u_{1}(x) \mapsto (2R)^{\frac{p + 1}{p - 1}} u_{1}(2R x).
\end{equation}
By small data arguments, $(\ref{4.1})$ implies that
\begin{equation}\label{4.3}
\| u \|_{L_{t,x}^{2(p - 1)}([0, \infty) \times \{ x : |x| \geq \frac{1}{2} + t \})} \lesssim \epsilon.
\end{equation}
Translating the initial data in time from $t = 0$ to $t = 1$,
\begin{equation}\label{4.4}
\| u \|_{L_{t,x}^{2(p - 1)}([1, \infty) \times \{ x : |x| \geq t - \frac{1}{2} \})} \lesssim \epsilon.
\end{equation}

As in \cite{D}, the proof of
\begin{equation}\label{4.5}
\| u \|_{L_{t,x}^{2(p - 1)}([1, \infty) \times \{ x : |x| \leq t - \frac{1}{2} \})} < \infty,
\end{equation}
will make use of the hyperbolic change of coordinates,
\begin{equation}\label{4.6}
\tilde{u}(\tau, s) = \frac{e^{\tau} \sinh s}{s} u(e^{\tau} \cosh s, e^{\tau} \sinh s).
\end{equation}
If $u$ solves $(\ref{1.1})$, then $\tilde{u}(\tau, s)$ solves
\begin{equation}\label{4.7}
(\partial_{\tau \tau} - \partial_{ss} - \frac{2}{s} \partial_{s}) \tilde{u}(\tau, s) + e^{-(p - 3)\tau} (\frac{s}{\sinh s})^{p - 1} |\tilde{u}(\tau, s)|^{p - 1} \tilde{u}(\tau, s) = 0.
\end{equation}

The hyperbolic energy is given by
\begin{equation}\label{4.8}
\aligned
E(\tilde{u}) = \frac{1}{2} \int (\partial_{s} \tilde{u}(\tau, s))^{2} s^{2} ds + \frac{1}{2} \int (\partial_{\tau} \tilde{u}(\tau, s))^{2} s^{2} ds \\ + \frac{1}{p + 1} \int e^{-(p - 3) \tau} (\frac{s}{\sinh s})^{p - 1} |\tilde{u}(\tau, s)|^{p + 1} s^{2} ds.
\endaligned
\end{equation}
By direct computation,
\begin{equation}\label{4.9}
\frac{d}{d\tau} E(\tilde{u})(\tau) = -\frac{p - 3}{p + 1} \int e^{-(p - 3) \tau} (\frac{s}{\sinh s})^{p - 1} |\tilde{u}(\tau, s)|^{p + 1} s^{2} ds \leq 0,
\end{equation}
which implies that the energy of $u$ is non-increasing.\vspace{5mm}

We also have a Morawetz estimate.
\begin{theorem}\label{t4.1}
If $\tilde{u}$ solves $(\ref{4.7})$ on any interval $I = [0, T]$, then
\begin{equation}\label{4.10}
\int_{I} \int e^{-(p - 3) \tau} (\frac{s}{\sinh s})^{p - 1} (\frac{\cosh s}{\sinh s}) |\tilde{u}(\tau, s)|^{p + 1} s^{2} ds d\tau \lesssim E(\tilde{u}(0)).
\end{equation}
\end{theorem}
\emph{Proof:} Again use the Morawetz potential in $(\ref{3.3})$,
\begin{equation}\label{4.11}
M(\tau) = \int \tilde{u}_{s}(s, \tau) \tilde{u}_{\tau}(s, \tau) s^{2} ds + \int \tilde{u}_{\tau}(s, \tau) \tilde{u}(s, \tau) s ds.
\end{equation}
Then by direct computation,
\begin{equation}\label{4.12}
\frac{d}{d\tau} M(\tau) = -\frac{1}{2} \tilde{u}(0, \tau)^{2} - \frac{p - 1}{p + 1} \int (\frac{\cosh s}{\sinh s}) (\frac{s}{\sinh s})^{p - 1} |\tilde{u}(s, \tau)|^{p + 1} ds.
\end{equation}
Then by $(\ref{4.9})$ and the fundamental theorem of calculus, the proof is complete. $\Box$

Previously, in \cite{D}, for the cubic wave equation, the initial data was split into a $(\tilde{v}_{0}, \tilde{v}_{1}) \in \dot{H}^{1} \times L^{2}$ component and a $(\tilde{w}_{0}, \tilde{w}_{1}) \in \dot{H}^{1/2} \times \dot{H}^{-1/2}$. Here, it would be nice if we could do something similar, only with $\dot{H}^{1/2}$ replaced by $\dot{H}^{s_{c}}$. However, this is not possible due to the fact that the hyperbolic energy scales like the $\dot{H}^{1/2}$ norm, and is not invariant under the general scaling $(\ref{1.2})$. Instead, what we will do is place $(\tilde{v}_{0}, \tilde{v}_{1}) \in \dot{H}^{1} \times L^{2}$, but $(\tilde{w}_{0}, \tilde{w}_{1})$ will merely lie in a Sobolev space after multiplying by exponential weights. The weights in the energy $(\ref{4.8})$ will then be used in conjunction with the weights for the Sobolev space to bound the growth of the energy of $\tilde{v}$.\vspace{5mm}

To calculate,
\begin{equation}\label{4.14.1}
\tilde{u}(\tau, s)|_{\tau = 0} = \frac{e^{\tau} \sinh s}{s} u(e^{\tau} \cosh s, e^{\tau} \sinh s)|_{\tau = 0},
\end{equation}
use Duhamel's principle,
\begin{equation}\label{4.15}
u(t) = S(t - 1)(u_{0}, u_{1}) + \int_{1}^{t} S(t - s)(0, |u|^{p - 1} u) ds.
\end{equation}
First consider the contribution of $S(t)(u_{0}, u_{1})$ with $u_{1} = 0$. In that case,
\begin{equation}\label{4.16}
\aligned
s \cdot S(t - 1)(u_{0}, 0)(e^{\tau} \cosh s, e^{\tau} \sinh s) \\ = \frac{1}{2} [u_{0}(e^{\tau + s} - 1) \cdot (e^{\tau + s} - 1) + u_{0}(1 - e^{\tau - s}) \cdot (1 - e^{\tau - s})].
\endaligned
\end{equation}
Again take $\phi \in C_{0}^{\infty}(\mathbf{R}^{3})$, only this time $\phi = 1$ on $|x| \leq 1$ and $\phi$ is supported on $|x| \leq 2$. Let $n$ be an integer satisfying $2^{n} \sim 2R$. By direct computation,
\begin{equation}\label{4.17}
\| \partial_{s} [\phi(e^{\tau + s} - 1) (P_{\leq n} u_{0})(e^{\tau + s} - 1) \cdot (e^{\tau + s} - 1)]|_{\tau = 0} \|_{L^{2}([0, \infty))} \lesssim 2^{n (1 - s_{c})} \| u_{0} \|_{\dot{H}^{s_{c}}(\mathbf{R}^{3})},
\end{equation}
and
\begin{equation}\label{4.18}
\|  [\phi(e^{\tau + s} - 1) (P_{\leq n} u_{0})(e^{\tau + s} - 1) \cdot (\frac{e^{\tau + s} - 1}{s})]|_{\tau = 0} \|_{L^{2}([0, \infty))} \lesssim 2^{n (1 - s_{c})} \| u_{0} \|_{\dot{H}^{s_{c}}(\mathbf{R}^{3})}.
\end{equation}
Meanwhile, by $(\ref{2.1.1})$,
\begin{equation}\label{4.19}
\| \phi(e^{\tau + s} - 1) (P_{\geq n} u_{0})(e^{\tau + s} - 1) \cdot (\frac{e^{\tau + s} - 1}{s})|_{\tau = 0} \|_{\dot{H}^{s_{c}}(\mathbf{R}^{3})} \lesssim \epsilon,
\end{equation}

Next, take a partition of unity,
\begin{equation}\label{4.21}
1 = \sum_{k \geq 0} \chi(s - k),
\end{equation}
where $\chi \in C_{0}^{\infty}(\mathbf{R})$, and $\chi$ is supported on $-1 \leq s \leq 1$. Then by direct computation,
\begin{equation}\label{4.22}
\aligned
\| \partial_{s} [\sum_{k \geq 0} \chi(s - k) \phi(1 - e^{\tau - s}) (P_{\leq n + \frac{k}{\ln(2)}} u_{0})(1 - e^{\tau - s}) \cdot (1 - e^{\tau - s})]|_{\tau = 0} \|_{L^{2}([0, \infty))} \\
\lesssim 2^{n(1 - s_{c})} \sum_{k \geq 0} \| u_{0} \|_{\dot{H}^{s_{c}}} e^{-k/2} e^{k (1 - s_{c})} + \sum_{k \geq 0} \| P_{n + \frac{k}{\ln(2)}} u_{0} \|_{\dot{H}^{1/2}} \lesssim 2^{n(1 - s_{c})} \| u_{0} \|_{\dot{H}^{s_{c}}},
\endaligned
\end{equation}
and
\begin{equation}\label{4.23}
\aligned
\| \sum_{k \geq 0} \chi(s - k) \phi(1 - e^{\tau - s}) (P_{\leq n + \frac{k}{\ln(2)}} u_{0})(1 - e^{\tau - s}) \cdot (\frac{1 - e^{\tau - s}}{s})|_{\tau = 0} \|_{L^{2}([0, \infty))} \\ \lesssim 2^{n(1 - s_{c})} (\sum_{k \geq 0} \frac{1}{1 + k^{2}})^{1/2} \| u_{0} \|_{\dot{H}^{s_{c}}(\mathbf{R}^{3})} \lesssim 2^{n(1 - s_{c})} \| u_{0} \|_{\dot{H}^{s_{c}}(\mathbf{R}^{3})}.
\endaligned
\end{equation}
Also by $(\ref{2.1.1})$,
\begin{equation}\label{4.24}
\| \sum_{k \geq 0} \chi(s - k) \phi(1 - e^{\tau - s}) (P_{> n + \frac{k}{\ln(2)}} u_{0})(1 - e^{\tau - s}) \cdot (\frac{1 - e^{\tau - s}}{s})|_{\tau = 0} \|_{\dot{H}^{s_{c}}(\mathbf{R}^{3})} \lesssim \epsilon.
\end{equation}
Finally, take
\begin{equation}\label{4.26}
[1 - \phi(e^{\tau + s} - 1)] u_{0}(e^{\tau + s} - 1) \cdot (\frac{e^{\tau + s} - 1}{s})|_{\tau = 0}.
\end{equation}
By a change of variables,
\begin{equation}\label{4.27}
\aligned
\| \chi(s - k) [1 - \phi(e^{\tau + s} - 1)] u_{0}(e^{\tau + s} - 1) \cdot (\frac{e^{\tau + s} - 1}{s})|_{\tau = 0} \|_{\dot{H}^{s_{c}}} \\ 
\lesssim e^{-k/2} (\int_{e^{k - 1}}^{e^{k + 1}} |u_{0}(r)|^{2} r^{2} dr)^{1/2} + e^{-k/2 + k \cdot s_{c}} (\int_{e^{k - 1}}^{e^{k + 1}} ||\nabla|^{s_{c}} u_{0}(r)|^{2} r^{2} dr)^{1/2}.
\endaligned
\end{equation}
Since $\partial_{\tau} f = \pm \partial_{s} f$ for the components of $(\ref{4.16})$, the same estimates also hold for $\partial_{\tau} \tilde{w}(\tau, s)|_{\tau = 0}$.\vspace{5mm}

We would like to use a Littlewood--Paley projection to split $(\ref{4.26})$ into a $\dot{H}^{1}$ component and a $\dot{H}^{s_{c}}$ component with appropriate bounds. The difficulty here is that the Littlewood--Paley projection is only known to have a rapidly decreasing weight, which does not commute well with an exponentially decreasing weight.

Instead, the estimate will rely on some frequency localized projection operators. Choose $\psi \in C_{0}^{\infty}(\mathbf{R}^{3})$ to be a radial, decreasing function supported on $|x| \leq \frac{1}{2}$, and such that $\int \psi(x) dx = 1$. Then define the Fourier multipliers
\begin{equation}\label{4.28}
\tilde{P}_{0} f(x) = \int \psi(x - y) f(y) dy,
\end{equation}
and for $j \geq 1$,
\begin{equation}\label{4.29}
\tilde{P}_{j} f(x) = 2^{3j} \int \psi(2^{j}(x - y)) f(y) dy - 2^{3(j - 1)} \int \psi(2^{j - 1}(x - y)) f(y) dy.
\end{equation}
Clearly,
\begin{equation}\label{4.29.1}
f = \sum_{j \geq 0} \tilde{P}_{j} f.
\end{equation}

Then ignoring the contribution of $(0, u_{1})$ and $|u|^{p - 1} u$ for a moment, let
\begin{equation}\label{4.29.2}
\tilde{v}_{0} = \sum_{k \geq 0} \sum_{j \leq n - \frac{k}{\ln(2)} \cdot \frac{s_{c} - \frac{1}{2}}{1 - s_{c}}} \tilde{P}_{j} \chi(s - k) \tilde{u}_{0}, \qquad \tilde{v}_{1} = \sum_{k \geq 0} \sum_{j \leq n - \frac{k}{\ln(2)} \cdot \frac{s_{c} - \frac{1}{2}}{1 - s_{c}}} \tilde{P}_{j} \chi(s - k) \tilde{u}_{1}.
\end{equation}
By $(\ref{4.22})$--$(\ref{4.27})$,
\begin{equation}\label{4.29.3}
\| \tilde{v}_{0} \|_{\dot{H}^{1}} + \| \tilde{v}_{1} \|_{L^{2}} \lesssim 2^{n(1 - s_{c})} \| u_{0} \|_{\dot{H}^{s_{c}}},
\end{equation}
and
\begin{equation}\label{4.29.4}
\sum_{k \geq 0} 2^{k(1 - 2s_{c})} \| \chi(s - k) \tilde{w}_{0} \|_{\dot{H}^{s_{c}}}^{2} + \sum_{k \geq 0} 2^{k(1 - 2 s_{c})} \| \chi(s - k) \tilde{w}_{1} \|_{\dot{H}^{s_{c} - 1}}^{2} \lesssim \epsilon^{2}.
\end{equation}

Turning to estimating the contribution of $S(t)(0, u_{1})$, split
\begin{equation}\label{4.29.5}
u_{1} = \phi(x) P_{\leq n} u_{1} + [u_{1} - \phi(x) P_{\leq n} u_{1}].
\end{equation}
By direct calculation,
\begin{equation}\label{4.29.6}
\| \partial_{\tau, s} \int_{1 - e^{-s}}^{e^{s} - 1} \phi(r) P_{\leq n} u_{1}(r) r dr \|_{L^{2}([0, \infty))} \lesssim 2^{n(1 - s_{c})} \| u_{1} \|_{\dot{H}^{s_{c} - 1}},
\end{equation}
and by H{\"o}lder's inequality,
\begin{equation}\label{4.29.7}
\| \frac{1}{s} \int_{1 - e^{-s}}^{e^{s} - 1} \phi(r) P_{\leq n} u_{1}(r) r dr \|_{L^{2}([0, \infty))} \lesssim 2^{n(1 - s_{c})} \| u_{1} \|_{\dot{H}^{s_{c} - 1}}.
\end{equation}

To handle the remainder in $(\ref{4.29.5})$, as in \cite{D}, observe that
\begin{equation}\label{4.30}
\frac{\sin(t \sqrt{-\Delta})}{\sqrt{-\Delta}} g = \partial_{t} (\frac{\cos(t \sqrt{-\Delta})}{\Delta} g).
\end{equation}
Plugging in the formula for a solution to the wave equation when $r > t$, let $w(t,r) = \cos(t \sqrt{-\Delta}) f$, where $f = \frac{g}{\Delta}$. Then,
\begin{equation}\label{4.31}
\aligned
\partial_{t} (w(t,r)) = \frac{1}{2r} \partial_{t} (f(t + r)(t + r) + f(r - t) (r - t)) \\ = \frac{1}{2r} [f(t + r) + f'(t + r)(t + r) - f(r - t) - f'(r - t) (r - t)].
\endaligned
\end{equation}
Since $f \in \dot{H}^{s_{c} + 1}(\mathbf{R}^{3})$, the contribution of
\begin{equation}\label{4.32}
f'(e^{\tau  + s} - 1) \cdot (e^{\tau + s} - 1)|_{\tau = 0}, \qquad f'(1 - e^{\tau - s}) \cdot (1 - e^{\tau - s})|_{\tau = 0}
\end{equation}
may be handled in a manner similar to the contribution of the terms arising from $S(t)(u_{0}, 0)$.

Meanwhile, by a change of variables, for $k \geq 1$,
\begin{equation}\label{4.33}
\aligned
\int (\chi(s - k) f'(e^{s} - 1) \cdot e^{s})^{2} ds \lesssim e^{2(s_{c} - \frac{1}{2}) k} \int_{e^{k - 1}}^{e^{k + 1}} |f'(r)|^{2} r^{2(1 - s_{c})} dr.
\endaligned
\end{equation}
and
\begin{equation}\label{4.34}
\int (\chi(s - k) f'(1 - e^{-s}) \cdot e^{-s})^{2} ds \lesssim e^{-k} (\int_{1 - e^{-k + 1}}^{1 - e^{-k - 1}} |f'(r)|^{2} dr \lesssim e^{-2k} \| f \|_{\dot{H}^{1 + s_{c}}(\mathbf{R}^{3})}^{2}.
\end{equation}
By an identical calculation,
\begin{equation}\label{4.35}
\int (\chi(s - k) \partial_{\tau} f(e^{s + \tau} - 1)|_{\tau = 0})^{2} ds = \int (\chi(s - k) f'(e^{s} - 1) \cdot e^{s})^{2} ds.
\end{equation}
and
\begin{equation}\label{4.36}
\int (\chi(s - k) \partial_{\tau} f(1 - e^{\tau - s})|_{\tau = 0})^{2} ds = \int (\chi(s - k) f'(1 - e^{-s}) \cdot e^{-s})^{2} ds \lesssim e^{-2k} \| f \|_{\dot{H}^{1 + s_{c}}(\mathbf{R}^{3})}^{2}.
\end{equation}
Also, by the fundamental theorem of calculus, for $s_{0} \sim 1$,
\begin{equation}\label{4.37}
s_{0} [f(e^{s_{0}} - 1) - f(1 - e^{-s_{0}})]^{2} = s_{0} [\int_{1 - e^{-s_{0}}}^{e^{s_{0}} - 1} f'(r) dr]^{2} \lesssim \int |f'(r)|^{2} r^{2(1 - s_{c})} dr \lesssim \| f \|_{\dot{H}^{s_{c}}}^{2}.
\end{equation}

Finally, consider
\begin{equation}\label{4.38}
f(e^{\tau + s} - 1) - f(1 - e^{\tau - s}),
\end{equation}
when $s < 1$. By direct computation,
\begin{equation}\label{4.39}
\partial_{\tau} [f(e^{\tau + s} - 1) - f(1 - e^{\tau - s})]|_{\tau = 0} = f'(e^{s} - 1) \cdot e^{s} + f'(1 - e^{-s}) \cdot e^{-s}.
\end{equation}
Then for $g \in \dot{H}^{1 - s_{c}}$, by Hardy's inequality,
\begin{equation}\label{4.40}
\int f'(e^{s} - 1) \cdot e^{s} \cdot g(s) s ds + \int f'(1 - e^{-s}) \cdot e^{-s} \cdot g(s) s ds \lesssim \| f \|_{\dot{H}^{1 + s_{c}}} \| g \|_{\dot{H}^{1 - s_{c}}}.
\end{equation}
Also, by the fundamental theorem of calculus,
\begin{equation}\label{4.41}
\aligned
f(e^{s} - 1) - f(1 - e^{-s}) = \int_{s - \frac{s^{2}}{2} + \frac{s^{3}}{3!} - ...}^{s + \frac{s^{2}}{2} + \frac{s^{3}}{3!} + ...} f'(r) dr \\
= \int_{0}^{1} f'(s + \theta (\frac{s^{2}}{2} + \frac{s^{3}}{3!} + ...)) \cdot (\frac{s^{2}}{2} + \frac{s^{3}}{3!} + ...) d\theta \\ + \int_{-1}^{0} f'(s + \theta(\frac{s^{2}}{2} - \frac{s^{3}}{3!} + ...) \cdot (\frac{s^{2}}{2} + \frac{s^{3}}{3!} + ...) d\theta.
\endaligned
\end{equation}
Therefore, since $\chi(s)$ is supported on $s \leq 1$,
\begin{equation}\label{4.42}
\| f(e^{s} - 1) - f(1 - e^{-s}) \|_{\dot{H}^{s_{c}}} \lesssim \| f \|_{\dot{H}^{s_{c}}}.
\end{equation}
Then take the contribution to $\tilde{v}_{0}$ to be
\begin{equation}\label{4.43}
\tilde{v}_{0} = \sum_{1 \leq k \leq \frac{n}{\ln(2)} \frac{1 - s_{c}}{s_{c} - \frac{1}{2}} } \chi(s - k) f(e^{\tau + s} - 1)|_{\tau = 0} + \sum_{k \geq 1} \chi(s - k) f(1 - e^{\tau - s})|_{\tau = 0},
\end{equation}
and the contribution to $\tilde{v}_{1}$ to be
\begin{equation}\label{4.43.1}
\tilde{v}_{1} = \sum_{1 \leq k \leq \frac{n}{\ln(2)} \frac{1 - s_{c}}{s_{c} - \frac{1}{2}}} \chi(s - k) \partial_{\tau} f(e^{\tau + s} - 1)|_{\tau = 0} + \sum_{k \geq 1} \chi(s - k) \partial_{\tau} f(1 - e^{\tau - s})|_{\tau = 0}.
\end{equation}

Now take the Duhamel term $u_{nl}$. Because the curve $t^{2} - r^{2} = 1$ has slope $\frac{dr}{dt} > 1$ everywhere,
\begin{equation}\label{4.44}
s \tilde{u}_{nl}(\tau, s)|_{\tau = 0} = \int_{1}^{e^{\tau} \cosh s} \int_{e^{\tau} \sinh s - e^{\tau} \cosh s + t}^{e^{\tau} \sinh s + e^{\tau} \cosh s - t} r |u|^{p - 1} u(t,r) dr dt.
\end{equation}
By direct computation,
\begin{equation}\label{4.45}
\aligned
\int_{0}^{k} (\partial_{s, \tau}(s \tilde{u}_{nl})|_{\tau = 0})^{2} ds \lesssim \int_{0}^{k} e^{2s} (\int_{1}^{\cosh s} (e^{s} - t) |u|^{p - 1} u(t, e^{s} - t) dt)^{2} ds \\ + \int_{0}^{k} e^{-2s} (\int_{1}^{\cosh s} (t - e^{-s}) |u|^{p - 1} u(t, t - e^{-s}) dt)^{2} ds.
\endaligned
\end{equation}
By H{\"o}lder's inequality, since $e^{s} - \cosh s \sim e^{s}$,
\begin{equation}\label{4.46}
\aligned
\int_{0}^{k} e^{2s} (\int_{1}^{\cosh s} (e^{s} - t) |u|^{p - 1} u(t, e^{s} - t) dt)^{2} ds \\ \lesssim \int_{0}^{k} \int_{1}^{\cosh s} e^{3s} (e^{s} - t)^{2} |u|^{2p}(t, e^{s} - t) dt ds \\ 
\lesssim \int_{0}^{e^{k}} \int_{t^{2} - r^{2} \leq 1} |u|^{2p}(t, r) r^{4} dt dr \lesssim e^{2(s_{c} - \frac{1}{2})k}.
\endaligned
\end{equation}
Also by a change of variables and H{\"o}lder's inequality, since $(t - e^{-s}) \gtrsim 1$ for $s \geq 1$ and $t \geq 1$,
\begin{equation}\label{4.46}
\aligned
\int_{1}^{\infty} e^{-2s} (\int_{1}^{\cosh s} (t - e^{-s}) |u|^{p - 1} u(t, t - e^{-s}) dt)^{2} ds \\ \lesssim \int_{1}^{\infty} \int_{1}^{\cosh s} e^{-s} (t - e^{-s})^{2} |u|^{2p}(t, e^{s} - t) dt ds \\ 
\lesssim \int_{2}^{\infty} \int_{t^{2} - r^{2} \leq 1} |u|^{2p}(t, r) r^{2} dt dr < \infty.
\endaligned
\end{equation}
Also, by the radial Sobolev embedding theorem and Young's inequality,
\begin{equation}\label{4.47}
\aligned
\int_{0}^{1} e^{-2s} (\int_{1}^{\cosh s} (t - e^{-s}) |u|^{p - 1} u(t, t - e^{-s}) dt)^{2} ds \\ \lesssim \int_{1}^{3} (\int_{t^{2} - r^{2} \leq 1} u(t, r)^{2p} r^{2} dr)^{1/2} dt \lesssim \int_{1}^{3} \frac{1}{(t - 1)^{-1 + \frac{s_{c}}{2}}} dt < \infty.
\endaligned
\end{equation}
This takes care of the nonlinear Duhamel piece.\vspace{5mm}

Therefore, we have finally proved:
\begin{theorem}
There exists a decomposition
\begin{equation}\label{4.48}
\tilde{u}_{0} = \tilde{v}_{0} + \tilde{w}_{0}, \qquad \tilde{u}_{1} = \tilde{v}_{1} + \tilde{w}_{1},
\end{equation}
satisfying
\begin{equation}\label{4.49}
\| \tilde{v}_{0} \|_{\dot{H}^{1}} + \| \tilde{v}_{1} \|_{L^{2}} \lesssim 2^{n(1 - s_{c})},
\end{equation}
and
\begin{equation}\label{4.50}
\sum_{k \geq 0} e^{(-2s_{c} + 1)k} \| \chi(s - k) \tilde{w}_{0} \|_{\dot{H}^{s_{c}} \cup \dot{H}^{1}}^{2} + \sum_{k \geq 0} e^{(-2 s_{c} + 1) k} \| \chi(s - k) \tilde{w}_{1} \|_{\dot{H}^{s_{c} - 1} \cup L^{2}}^{2} \lesssim \epsilon^{2}.
\end{equation}

\end{theorem}

\section{Scattering : Virial identities}
Now we are ready to prove scattering.
\begin{theorem}\label{t5.1}
For any radial $(u_{0}, u_{1})$, the global solution to $(\ref{1.1})$ scatters both forward and backward in time.
\end{theorem}
\emph{Proof:} Modifying $(\ref{3.13})$ and $(\ref{3.14})$, split $\tilde{u} = \tilde{v} + \tilde{w}$, where $\tilde{w}$ solves
\begin{equation}\label{4.13}
\aligned
\partial_{\tau \tau} \tilde{w} - \partial_{ss} \tilde{w} - \frac{2}{s} \partial_{s} \tilde{w} + e^{-(p - 3) \tau} (\frac{s}{\sinh s})^{p - 1} \cdot \tilde{P}_{\geq n - \frac{s_{c} - \frac{1}{2}}{1 - s_{c}} \frac{\tau}{\ln(2)}} |\tilde{w}|^{p - 1} \tilde{w} = 0, \\ \qquad w(0, s) = \tilde{w}_{0}(s), \qquad w_{\tau}(0,s) = \tilde{w}_{1},
\endaligned
\end{equation}
and $\tilde{v}$ solves
\begin{equation}\label{4.14}
\aligned
\partial_{\tau \tau} \tilde{v} - \partial_{ss} \tilde{v} - \frac{2}{s} \partial_{s} \tilde{v} + e^{-(p - 3) \tau} (\frac{s}{\sinh s})^{p - 1} [|\tilde{u}|^{p - 1} \tilde{u} - \tilde{P}_{\geq n - \frac{s_{c} - \frac{1}{2}}{1 - s_{c}} \frac{\tau}{\ln(2)}} |\tilde{w}|^{p - 1} \tilde{w}] = 0, \\
\tilde{v}(0, s) = \tilde{v}_{0}, \qquad \tilde{v}_{\tau}(0, s) = \tilde{v}_{1},
\endaligned
\end{equation}
where $\tilde{u}_{0} = \tilde{w}_{0} + \tilde{v}_{0}$ and $\tilde{u}_{1} = \tilde{w}_{1} + \tilde{v}_{1}$.\vspace{5mm}

$(\ref{4.13})$ may be shown to be scattering using small data arguments. Indeed, Strichartz estimates, finite propagation speed, $(\ref{4.50})$, and the fact that $(\frac{s}{\sinh s})$ is rapidly decreasing in $s$, imply that
\begin{equation}\label{5.7}
\aligned
 \| e^{-\frac{p - 3}{p - 1} \tau} (\frac{s}{\sinh s})^{p - 1} e^{-(p - 3) \tau} |\tilde{w}|^{p - 1} \tilde{w} \|_{L_{\tau}^{\frac{2}{1 + s_{c}}} L_{s}^{\frac{2}{2 - s_{c}}}} \\
\lesssim \| e^{-\frac{p - 3}{p - 1} \tau} (\frac{s}{\sinh s})^{\frac{p - 3}{p - 1}+} \tilde{w} \|_{L_{\tau, s}^{2(p - 1)}}^{p - 1} \| e^{-\frac{p - 3}{p - 1} \tau} (\frac{s}{\sinh s})^{\frac{p - 3}{p - 1}+} \tilde{w} \|_{L_{\tau}^{\frac{2}{s_{c}}} L_{s}^{\frac{2}{1 - s_{c}}}} \\
\lesssim \epsilon^{p} + \| e^{-\frac{p - 3}{p - 1} \tau} (\frac{s}{\sinh s})^{p - 1} e^{-(p - 3) \tau} |\tilde{w}|^{p - 1} \tilde{w} \|_{L_{\tau}^{\frac{2}{1 + s_{c}}} L_{s}^{\frac{2}{2 - s_{c}}}}^{p}.
\endaligned
\end{equation}
These estimates also commute well with the $\tilde{P}_{j}$ operators.
The same calculation also shows
\begin{equation}\label{5.7.1}
\| e^{-\frac{p - 3}{p - 1} \tau} e^{-(p - 3) \tau} (\frac{s}{\sinh s})^{p - 1} |\tilde{w}|^{p - 1} \tilde{w} \|_{L_{\tau}^{1} L_{s}^{\frac{6}{5 - 2 s_{c}}}} \lesssim \epsilon^{p}.
\end{equation}

Now, define the modified energy
\begin{equation}\label{5.1}
\mathcal E(\tau) = E(\tau) + c M(\tau) + \int \sum_{j \geq 0} [|\tilde{v}|^{p - 1} \tilde{v} - |\tilde{P}_{\leq j} \tilde{v}|^{p - 1} (\tilde{P}_{j} \tilde{v})] \cdot \tilde{P}_{j} \tilde{w} ds,
\end{equation}
where $c > 0$ is a small constant,
\begin{equation}\label{5.2}
\aligned
E(\tau) = \frac{1}{2} \int \tilde{v}_{s}(s, \tau)^{2} s^{2} ds + \frac{1}{2} \int \tilde{v}_{\tau}(s, \tau)^{2} s^{2} ds \\ + \frac{1}{p + 1} \int e^{-(p - 3)\tau} (\frac{s}{\sinh s})^{p - 1} |\tilde{v}(s, \tau)|^{p + 1} s^{2} ds,
\endaligned
\end{equation}
and
\begin{equation}\label{5.3}
M(\tau) = \int \tilde{v}_{s}(s, \tau) \tilde{v}_{\tau}(s, \tau) s^{2} ds + \int \tilde{v}(s, \tau) \tilde{v}_{\tau}(s, \tau) s^{2} ds.
\end{equation}
Then by $(\ref{4.9})$ and $(\ref{4.12})$,
\begin{equation}\label{5.4}
\aligned
\frac{d}{d\tau} \mathcal E(\tau) = -\frac{c}{2} \tilde{v}(\tau, 0)^{2} - c \frac{p - 1}{p + 1} \int (\frac{\cosh s}{\sinh s}) e^{-(p - 3) \tau} (\frac{s}{\sinh s})^{p - 1} |\tilde{v}(s, \tau)|^{p + 1} s^{2} ds \\
- \frac{p - 3}{p + 1} \int (\frac{s}{\sinh s})^{p - 1} e^{-(p - 3) \tau} |\tilde{v}(s, \tau)|^{p + 1} s^{2} ds \\
+ \frac{d}{d \tau}  \int \sum_{j \geq 0} [|\tilde{v}|^{p - 1} \tilde{v} - |\tilde{P}_{\leq j} \tilde{v}|^{p - 1} (\tilde{P}_{j} \tilde{v})] \cdot \tilde{P}_{j} \tilde{w} s^{2} ds \\ - \int e^{-(p - 3) \tau} (\frac{s}{\sinh s})^{p - 1} [|\tilde{u}|^{p - 1} \tilde{u} - |\tilde{v}|^{p - 1} \tilde{v} - \tilde{P}_{\geq n - \frac{s_{c} - \frac{1}{2}}{1 - s_{c}} \frac{\tau}{\ln(2)}}  |\tilde{w}|^{p - 1} \tilde{w}] \tilde{v}_{\tau} s^{2} ds \\
- \int e^{-(p - 3)\tau} (\frac{s}{\sinh s})^{p - 1} [|\tilde{u}|^{p - 1} \tilde{u} - |\tilde{v}|^{p - 1} \tilde{v} - \tilde{P}_{\geq n - \frac{s_{c} - \frac{1}{2}}{1 - s_{c}} \frac{\tau}{\ln(2)}} |\tilde{w}|^{p - 1} \tilde{w}] \tilde{v}_{s} s^{2} ds \\
- \int e^{-(p - 3) \tau} (\frac{s}{\sinh s})^{p - 1} [|\tilde{u}|^{p - 1} \tilde{u} - |\tilde{v}|^{p - 1} \tilde{v} - \tilde{P}_{\geq n - \frac{s_{c} - \frac{1}{2}}{1 - s_{c}} \frac{\tau}{\ln(2)}}  |\tilde{w}|^{p - 1} \tilde{w}] \tilde{v} s ds.
\endaligned
\end{equation}
First consider the contribution of $P_{\leq n - \frac{s_{c} - \frac{1}{2}}{1 - s_{c}} \frac{\tau}{\ln(2)}} |\tilde{w}|^{p - 1} \tilde{w}$. By the Sobolev embedding theorem,
\begin{equation}\label{5.4.1}
\aligned
\| e^{-(p - 3) \tau} (\frac{s}{\sinh s}) P_{\leq n - \frac{s_{c} - \frac{1}{2}}{1 - s_{c}}} (|\tilde{w}|^{p - 1} \tilde{w}) \cdot \tilde{v}_{s, \tau} \|_{L^{1}} \\
\lesssim 2^{n(1 - s_{c})} E(\tilde{v})^{1/2} \| e^{-(p - 3) \tau} (\frac{s}{\sinh s})^{p - 2} |\tilde{w}|^{p - 1} \|_{L_{\tau, s}^{2}} \| e^{-\frac{p - 3}{p - 1} \tau} (\frac{s}{\sinh s}) \tilde{w} \|_{L^{\frac{3}{1 - s_{c}}}}.
\endaligned
\end{equation}

Next, as in $(\ref{3.21})$,
\begin{equation}\label{5.5}
\aligned
\int e^{-(p - 3) \tau} (\frac{s}{\sinh s})^{p - 1} [|\tilde{u}|^{p - 1} \tilde{u} - |\tilde{v}|^{p - 1} \tilde{v} - |\tilde{w}|^{p - 1} \tilde{w}] \tilde{v} s ds \\
\lesssim e^{-(p - 3) \tau} (\frac{s}{\sinh s})^{p - 1} [|\tilde{v}|^{p} |\tilde{w}| + |\tilde{w}|^{p - 1} |\tilde{v}|^{2}] s ds \\
\lesssim (\int e^{-(p - 3) \tau} (\frac{s}{\sinh s})^{p - 1} |\tilde{v}|^{p + 1} (\frac{\cosh s}{\sinh s}) s^{2} ds)^{\frac{p - 2}{p - 1}} \| \frac{1}{s^{1/2}} \tilde{v} \|_{L^{3}}^{\frac{2}{p - 1}} \| e^{-\frac{p - 3}{p - 1} \tau} (\frac{s}{\sinh s}) \tilde{w} \|_{L^{3(p - 1)}} \\
+ \| e^{-\frac{p - 3}{p - 1} \tau} (\frac{s}{\sinh s}) \tilde{w} \|_{L^{3(p - 1)}}^{p - 1} \| \frac{1}{s} \tilde{v} \|_{L^{2}} \| \tilde{v} \|_{L^{6}} \\
\lesssim \delta (\int e^{-(p - 3) \tau} (\frac{s}{\sinh s})^{p - 1} |\tilde{v}|^{p + 1} s^{2} ds) + \frac{1}{\delta} E(\tau) \| e^{-\frac{p - 3}{p - 1} \tau} (\frac{s}{\sinh s}) \tilde{w} \|_{L^{3(p - 1)}}^{p - 1}.
\endaligned
\end{equation}

Next, following $(\ref{3.23})$ and $(\ref{3.24})$,
\begin{equation}\label{5.8}
\aligned
\int e^{-(p - 3) \tau} (\frac{s}{\sinh s})^{p - 1} |\tilde{v}|^{p - 2} |\partial_{s,\tau} \tilde{v}| |\tilde{w}|^{2} s^{2} ds \\ \lesssim \| e^{-\frac{p - 3}{p - 1} \tau} (\frac{s}{\sinh s}) \tilde{w} \|_{L^{\frac{3}{1 - s_{c}}}}^{2} \| e^{-\frac{p - 3}{p + 1} \tau} (\frac{s}{\sinh s})^{\frac{p - 1}{p + 1}} \tilde{v} \|_{L^{p + 1}}^{p - 2 - \frac{5 - p}{p - 1}} \| \partial_{s, \tau} \tilde{v} \|_{L^{2}}^{1 + \frac{5 - p}{p - 1}} \\
\lesssim E(\tau)  \| e^{-\frac{p - 3}{p - 1} \tau} (\frac{s}{\sinh s}) \tilde{w} \|_{L^{\frac{3}{1 - s_{c}}}}^{2}.
\endaligned
\end{equation}
Also by Strichartz estimates, the weights $e^{-(p - 3) \tau} (\frac{s}{\sinh s})^{p - 1}$, and $(\ref{4.50})$,
\begin{equation}\label{5.9}
\int \| e^{-\frac{p - 3}{p - 1} \tau} (\frac{s}{\sinh s}) \tilde{w} \|_{L^{\frac{3}{1 - s_{c}}}}^{2} < \epsilon^{2}.
\end{equation}
Meanwhile,
\begin{equation}\label{5.10}
\aligned
\int e^{-(p - 3) \tau} (\frac{s}{\sinh s})^{p - 1} |\tilde{w}|^{p - 1} |\tilde{v}| |\partial_{s, \tau} \tilde{v}| s^{2} ds \\ \lesssim \| e^{-\frac{p - 3}{p - 1} \tau} (\frac{s}{\sinh s}) \tilde{w} \|_{L^{3(p - 1)}}^{p - 1} \| \partial_{s, \tau} \tilde{v} \|_{L^{2}}^{2} \\ \lesssim E(\tau) \| e^{-\frac{p - 3}{p - 1} \tau} (\frac{s}{\sinh s}) \tilde{w} \|_{L^{3(p - 1)}}^{p - 1}.
\endaligned
\end{equation}
By $(\ref{5.7})$ and $(\ref{5.9})$, these terms may be handled using Gronwall's inequality.

Now then, by the product rule,
\begin{equation}\label{5.11}
\aligned
\frac{d}{d\tau} \int \sum_{j} [|\tilde{v}|^{p - 1} \tilde{v} - |\tilde{P}_{\leq j} \tilde{v}|^{p - 1} \tilde{P}_{\leq j} \tilde{v}] e^{-(p - 3) \tau} (\frac{s}{\sinh s})^{p - 1} \tilde{P}_{j} \tilde{w} s^{2} ds \\ - p \int e^{-(p - 3) \tau} (\frac{s}{\sinh s})^{p - 1} |\tilde{v}|^{p - 1} \tilde{v}_{\tau} \tilde{w} s^{2} ds \\
= -\int \sum_{j} e^{-(p - 3) \tau} (\frac{s}{\sinh s})^{p - 1}  [|\tilde{v}|^{p - 1} \tilde{v} - |\tilde{P}_{\leq j} \tilde{v}|^{p - 1} \tilde{P}_{\leq j} \tilde{v}] \tilde{w}_{\tau} s^{2} ds \\
- \int \sum_{j} e^{-(p - 3) \tau} (\frac{s}{\sinh s})^{p - 1} \partial_{\tau} [|\tilde{P}_{\leq j} \tilde{v}|^{p - 1} \tilde{P}_{\leq j} \tilde{v}] \tilde{P}_{j} \tilde{w} s^{2} ds \\
- (p - 3) \int |\tilde{v}|^{p - 1} \tilde{v} \tilde{w} e^{-(p - 3) \tau} (\frac{s}{\sinh s})^{p - 1} s^{2} ds.
\endaligned
\end{equation}
By the Cauchy--Schwartz inequality,
\begin{equation}\label{5.12}
\aligned
- (p - 3) \int |\tilde{v}|^{p - 1} \tilde{v} \tilde{w} e^{-(p - 3) \tau} (\frac{s}{\sinh s})^{p - 1} s^{2} ds \lesssim \delta (\int |\tilde{v}|^{p + 1} e^{-(p - 3) \tau} (\frac{s}{\sinh s})^{p - 1} s^{2} ds) \\ + \frac{1}{\delta} \| e^{-\frac{p - 3}{p - 1} \tau} (\frac{s}{\sinh s}) \tilde{w} \|_{L^{3(p - 1)}}^{p - 1} \| \partial_{s, \tau} \tilde{v} \|_{L^{2}}^{2}.
\endaligned
\end{equation}


Next,
\begin{equation}\label{5.12.1}
\aligned
\int e^{-(p - 3) \tau} (\frac{s}{\sinh s})^{p - 1} \partial_{\tau} [|\tilde{P}_{\leq j} \tilde{v}|^{p - 1} \tilde{P}_{\leq j} \tilde{v}] \tilde{P}_{j} \tilde{w} s^{2} ds \\
\lesssim \| (\frac{\sinh s}{\cosh s})^{\frac{p - 3}{2(p + 1)}} \tilde{P}_{\leq j} \tilde{v}_{\tau} \|_{L^{\frac{p + 1}{2}}} (\int e^{-(p - 3) \tau} (\frac{s}{\sinh s})^{p - 1} (\frac{\cosh s}{\sinh s}) |\tilde{v}|^{p + 1} s^{2} ds)^{\frac{p - 1}{p + 1}} \\ \times \| (\frac{\sinh s}{\cosh s})^{1/2} e^{-\frac{2(p - 3)}{p + 1} \tau} (\frac{s}{\sinh s})^{\frac{2(p - 1)}{p + 1}} \tilde{P}_{j} \tilde{w}_{\tau} \|_{L^{\infty}}.
\endaligned
\end{equation}
By the radial Sobolev embedding theorem,
\begin{equation}\label{5.12.2}
\| (\frac{\sinh s}{\cosh s})^{\frac{p - 3}{2(p + 1)}} \tilde{P}_{\leq j} \tilde{v}_{\tau} \|_{L^{\frac{p + 1}{2}}} \lesssim 2^{\frac{j(p - 3)}{p + 1}} E(\tilde{v})^{1/2}.
\end{equation}

By direct computation,
\begin{equation}\label{5.16}
[|\tilde{v}|^{p - 1} \tilde{v} - |\tilde{P}_{\leq j} \tilde{v}|^{p - 1} \tilde{P}_{\leq j} \tilde{v}] = O(|\tilde{P}_{> j} \tilde{v}| (|\tilde{P}_{\leq j} \tilde{v}|^{p - 1} + |\tilde{P}_{> j} \tilde{v}|^{p - 1})).
\end{equation}
Also, by Bernstein's inequality and the radial Sobolev embedding theorem,
\begin{equation}\label{5.17}
\aligned
\int e^{-(p - 3) \tau} (\frac{s}{\sinh s})^{p - 1} [|\tilde{v}|^{p - 1} \tilde{v} - |\tilde{P}_{\leq j} \tilde{v}|^{p - 1} \tilde{P}_{\leq j} \tilde{v}] \tilde{P}_{j} \tilde{w}_{\tau} s^{2} ds \\
\lesssim (\int e^{-(p - 3) \tau} (\frac{s}{\sinh s})^{p - 1} (\frac{\cosh s}{\sinh s}) |\tilde{v}|^{p + 1} s^{2} ds)^{\frac{p - 1}{p + 1}} \\ \times \| (\frac{\sinh s}{\cosh s})^{\frac{p - 3}{2(p + 1)}} |\tilde{P}_{> j} \tilde{v}| \|_{L^{\frac{p + 1}{2}}} \| (\frac{\sinh s}{\cosh s})^{1/2} e^{-\frac{2(p - 3)}{p + 1} \tau} (\frac{s}{\sinh s})^{\frac{2(p - 1)}{p + 1}} \tilde{P}_{j} \tilde{w} \|_{L^{\infty}}.
\endaligned
\end{equation}
By the radial Sobolev embedding theorem and the definition of $\tilde{P}_{j}$,
\begin{equation}\label{5.18}
\| (\frac{\sinh s}{\cosh s})^{\frac{p - 3}{2(p + 1)}} |\tilde{P}_{> j} \tilde{v}| \|_{L^{\frac{p + 1}{2}}} \lesssim 2^{-\frac{4j}{p + 1}} E(\tilde{v})^{1/2}.
\end{equation}
Now then, summing up the contribution of the linear term to $\tilde{w}$,
\begin{equation}\label{5.19}
\sum_{j} \sum_{k \geq n \frac{1 - s_{c}}{s_{c} - \frac{1}{2}} - j \frac{1 - s_{c}}{s_{c} - \frac{1}{2}}} 2^{-k \frac{2(p - 3)}{p + 1} + k \frac{p - 3}{p - 1}} 2^{j(1 - s_{c} + \frac{p - 1}{p + 1})} 2^{\frac{-4j}{p + 1}} \lesssim \epsilon 2^{n(1 - s_{c}) \cdot \frac{3 - p}{p + 1}} \lesssim \epsilon E(0)^{\frac{3 - p}{2(p + 1)}}.
\end{equation}
Also, considering the contribution of the nonlinear term,
\begin{equation}\label{5.20}
\aligned
e^{-\frac{2(p - 3)}{p + 1} \tau} e^{\frac{p - 3}{p - 1} \tau} \|  e^{-\tau (p - 3)} (\frac{s}{\sinh s})^{p - 1} \tilde{P}_{\geq n - \frac{\tau}{\ln(2)} \frac{s_{c} - \frac{1}{2}}{1 - s_{c}}} |\tilde{w}|^{p - 1} \tilde{w} d\tau \|_{\dot{H}^{s_{c} - \frac{(3 - p)(1 - s_{c})}{(p + 1)}}}  \\\lesssim E(0)^{\frac{(3 - p)}{2(p + 1)}} \| e^{-\tau \frac{p - 3}{p - 1}} e^{-\tau (p - 3)} (\frac{s}{\sinh s})^{p - 1} |\tilde{w}|^{p - 1} \tilde{w} \|_{L^{\frac{6}{5 - s_{c}}}}.
\endaligned
\end{equation}
Then by $(\ref{5.7.1})$ and a Gronwall-type estimate, we have proved
\begin{equation}\label{5.21}
\int \int |\tilde{v}(\tau, s)|^{p + 1} (\frac{\cosh s}{\sinh s}) (\frac{s}{\sinh s})^{p - 1} e^{-(p - 3) \tau} s^{2} ds d\tau < \infty.
\end{equation}
By the radial Sobolev embedding theorem,
\begin{equation}\label{5.22}
(\frac{\sinh s}{\cosh s}) |\tilde{v}(\tau, s)|^{p - 3} \lesssim E(\tilde{v})^{\frac{p - 3}{2}} < \infty.
\end{equation}
Therefore, we have proved
\begin{equation}\label{5.23}
\int \int |v(e^{\tau} \cosh s, e^{\tau} \sinh s)|^{2(p - 1)} (e^{\tau} \sinh s)^{2} e^{2\tau} ds d\tau < \infty,
\end{equation}
which by a change of variables formula implies
\begin{equation}\label{5.24}
\int_{1}^{\infty} \int_{t^{2} - r^{2} \geq 1} |v(t, r)|^{2(p - 1)} r^{2} dr dt < \infty.
\end{equation}
Also, by $(\ref{5.7})$ and a change of variables,
\begin{equation}\label{5.25}
\aligned
\int \int e^{-2(p - 3) \tau} |\tilde{w}(e^{\tau} \cosh s, e^{\tau} \sinh s)|^{2(p - 1)} (\frac{s}{\sinh s})^{2(p - 2)} s^{2} ds d\tau \\
= \int \int e^{2 \tau} |w(e^{\tau} \cosh s, e^{\tau} \sinh s)|^{2(p - 1)} (e^{\tau} \sinh s)^{2} ds d\tau \\
= \int \int_{t^{2} - r^{2} \geq 1} |w(t, r)|^{2(p - 1)} r^{2} dr dt \lesssim \epsilon^{2(p - 1)}.
\endaligned
\end{equation}
Combining $(\ref{4.4})$ with $(\ref{5.24})$ and $(\ref{5.25})$ completes the proof of Theorem $\ref{t5.1}$. $\Box$

\section{Scattering}
As in \cite{D2} and \cite{D}, $(\ref{1.4.1})$ is proved by combining Zorn's lemma and a perturbative argument.

Let $(u_{0}^{n}, u_{1})$ be a uniformly bounded radial sequence,
\begin{equation}\label{6.1}
\| u_{0}^{n} \|_{\dot{H}^{s_{c}}(\mathbf{R}^{3})} + \| u_{1}^{n} \|_{\dot{H}^{s_{c} - 1}(\mathbf{R}^{3})} \leq A,
\end{equation}
and let $u^{n}$ be the solution to $(\ref{1.1})$ with initial data $(u_{0}^{n}, u_{1}^{n})$. By Zorn's lemma,   to prove $(\ref{1.4.1})$, it suffices to show that
\begin{equation}\label{6.2}
\| u^{n} \|_{L_{t,x}^{2(p - 1)}(\mathbf{R} \times \mathbf{R}^{3})}
\end{equation}
is uniformly bounded for any such sequence.

The proof of this fact makes use of the profile decomposition.
\begin{theorem}[Profile decomposition]\label{t6.1}
Suppose that there is a uniformly bounded, radially symmetric sequence
\begin{equation}\label{6.3}
\| u_{0}^{n} \|_{\dot{H}^{s_{c}}(\mathbf{R}^{3})} + \| u_{1}^{n} \|_{\dot{H}^{s_{c} - 1}(\mathbf{R}^{3})} \leq A < \infty.
\end{equation}
Then there exists a subsequence, also denoted $(u_{0}^{n}, u_{1}^{n}) \subset \dot{H}^{s_{c}} \times \dot{H}^{s_{c} - 1}$ such that for any $N < \infty$,
\begin{equation}\label{6.4}
S(t)(u_{0}^{n}, u_{1}^{n}) = \sum_{j = 1}^{N} \Gamma_{n}^{j} S(t)(\phi_{0}^{j}, \phi_{1}^{j}) + S(t)(R_{0, n}^{N}, R_{1,n}^{N}),
\end{equation}
with
\begin{equation}\label{6.5}
\lim_{N \rightarrow \infty} \limsup_{n \rightarrow \infty} \| S(t)(R_{0,n}^{N}, R_{1,n}^{N}) \|_{L_{t,x}^{q}(\mathbf{R} \times \mathbf{R}^{3})} = 0.
\end{equation}
$\Gamma_{n}^{j} = (\lambda_{n}^{j}, t_{n}^{j})$ belongs to the group $(0, \infty) \times \mathbf{R}$, which acts by
\begin{equation}\label{6.6}
\Gamma_{n}^{j} F(t,x) = \lambda_{n}^{j} F(\lambda_{n}^{j} (t - t_{n}^{j}), \lambda_{n}^{j} x).
\end{equation}
The $\Gamma_{n}^{j}$ are pairwise orthogonal, that is, for every $j \neq k$,
\begin{equation}\label{6.7}
\lim_{n \rightarrow \infty} \frac{\lambda_{n}^{j}}{\lambda_{n}^{k}} + \frac{\lambda_{n}^{k}}{\lambda_{n}^{j}} + (\lambda_{n}^{j})^{1/2} (\lambda_{n}^{k})^{1/2} |t_{n}^{j} - t_{n}^{k}| = \infty.
\end{equation}
Furthermore, for every $N \geq 1$,
\begin{equation}\label{6.8}
\aligned
\| (u_{0, n}, u_{1, n}) \|_{\dot{H}^{s_{c}} \times \dot{H}^{s_{c} - 1}}^{2} = \sum_{j = 1}^{N} \| (\phi_{0}^{j}, \phi_{0}^{k}) \|_{\dot{H}^{s_{c}} \times \dot{H}^{s_{c} - 1}}^{2} \\ + \| (R_{0, n}^{N}, R_{1, n}^{N}) \|_{\dot{H}^{s_{c}} \times \dot{H}^{s_{c} - 1}}^{2} + o_{n}(1).
\endaligned
\end{equation}
\end{theorem}

Theorem $\ref{t6.1}$ gives the profile decomposition
\begin{equation}\label{6.9}
S(t)(u_{0}^{n}, u_{1}^{n}) = \sum_{j = 1}^{N} S(t - t_{n}^{j}) (\lambda_{n}^{j} \phi_{0}^{j}(\lambda_{n}^{j} x), (\lambda_{n}^{j})^{2} \phi_{1}^{j}(\lambda_{n}^{j} x)) + S(t)(R_{0, n}^{N}, R_{1,n}^{N}),
\end{equation}
and moreover,
\begin{equation}\label{6.10}
S(\lambda_{n}^{j} t_{n}^{j})(\frac{1}{\lambda_{n}^{j}} u_{0}^{n}(\frac{x}{\lambda_{n}^{j}}), \frac{1}{(\lambda_{n}^{j})^{2}} u_{1}^{n}(\frac{x}{\lambda_{n}^{j}})) \rightharpoonup \phi_{0}^{j}(x),
\end{equation}
weakly in $\dot{H}^{s_{c}}(\mathbf{R}^{3})$, and
\begin{equation}\label{6.11}
\partial_{t}S(t + \lambda_{n}^{j} t_{n}^{j})(\frac{1}{\lambda_{n}^{j}} u_{0}^{n}(\frac{x}{\lambda_{n}^{j}}), \frac{1}{(\lambda_{n}^{j})^{2}} u_{1}^{n}(\frac{x}{\lambda_{n}^{j}}))|_{t = 0} \rightharpoonup \phi_{1}^{j}(x)
\end{equation}
weakly in $\dot{H}^{s_{c} - 1}(\mathbf{R}^{3})$. 

First consider the case that $\lambda_{n}^{j} t_{n}^{j}$ is uniformly bounded. In this case, after passing to a subsequence, $\lambda_{n}^{j} t_{n}^{j}$ converges to some $t^{j}$. Changing $(\phi_{0}^{j}, \phi_{1}^{j})$ to $S(-t^{j})(\phi_{0}^{j}, \phi_{1}^{j})$ and absorbing the error into $(R_{0, n}^{N}, R_{1, n}^{N})$,
\begin{equation}\label{6.12}
(\frac{1}{\lambda_{n}^{j}} u_{0}^{n}(\frac{x}{\lambda_{n}^{j}}), \frac{1}{(\lambda_{n}^{j})^{2}} u_{1}^{n}(\frac{x}{\lambda_{n}^{j}})) \rightharpoonup \phi_{0}^{j}(x),
\end{equation}
and
\begin{equation}\label{6.13}
\partial_{t}S(t)(\frac{1}{\lambda_{n}^{j}} u_{0}^{n}(\frac{x}{\lambda_{n}^{j}}), \frac{1}{(\lambda_{n}^{j})^{2}} u_{1}^{n}(\frac{x}{\lambda_{n}^{j}}))|_{t = 0} \rightharpoonup \phi_{1}^{j}(x).
\end{equation}
Then if $u^{j}$ is the solution to $(\ref{1.1})$ with initial data $(\phi_{0}^{j}, \phi_{1}^{j})$, then
\begin{equation}\label{6.14}
\| u^{j} \|_{L_{t,x}^{2(p - 1)}(\mathbf{R} \times \mathbf{R}^{3})} \leq M_{j} < \infty.
\end{equation}

Next, suppose that after passing to a subsequence, $\lambda_{n}^{j} t_{n}^{j} \nearrow +\infty$. Then a solution to $(\ref{1.1})$ approaches a translation in time of a solution to $(\ref{1.1})$ that scatters backward in time to $S(t)(\phi_{0}, \phi_{1})$, that is,
\begin{equation}\label{6.15}
\lim_{t \rightarrow -\infty} \| u - S(t)(\phi_{0}, \phi_{1}) \|_{\dot{H}^{s_{c}} \times \dot{H}^{s_{c} - 1}} = 0.
\end{equation}
Indeed, by Strichartz estimates, the dominated convergence theorem, and small data arguments, for some $T < \infty$ sufficiently large, $(\ref{1.1})$ has a solution $u$ on $(-\infty, -T]$ such that
\begin{equation}\label{6.16}
\| u \|_{L_{t,x}^{2(p - 1)} \cap L_{t}^{\frac{2}{s_{c}}} L_{x}^{\frac{2}{1 - s_{c}}}((-\infty, -T] \times \mathbf{R}^{3})} \lesssim \epsilon, \qquad (u(-T, x), u_{t}(-T, x)) = S(-T)(\phi_{0}, \phi_{1}),
\end{equation}
and by Strichartz estimates,
\begin{equation}\label{6.17}
\lim_{t \rightarrow +\infty} \| S(t)(u(-t), u_{t}(-t)) - (\phi_{0}, \phi_{1}) \|_{\dot{H}^{s_{c}} \times \dot{H}^{s_{c} - 1}} \lesssim \epsilon^{p}.
\end{equation}
Then by the inverse function theorem, there exists some $(u_{0}(-T), u_{1}(-T))$ such that $(\ref{1.1})$ has a solution that scatters backward in time to $S(t)(\phi_{0}, \phi_{1})$. Moreover, this solution must also scatter forward in time. Therefore,
\begin{equation}\label{6.18}
S(-t_{n}^{j})(\lambda_{n}^{j} \phi_{0}^{j}(\lambda_{n}^{j} x), (\lambda_{n}^{j})^{2} \phi_{1}^{j}(\lambda_{n}^{j} x))
\end{equation}
converges strongly to
\begin{equation}\label{6.19}
(\lambda_{n}^{j} u^{j}(-\lambda_{n}^{j} t_{n}^{j}, \lambda_{n}^{j} x), (\lambda_{n}^{j})^{2} u_{t}^{j}(-\lambda_{n}^{j} t_{n}^{j}, \lambda_{n}^{j} x))
\end{equation}
in $\dot{H}^{s_{c}} \times \dot{H}^{s_{c} - 1}$, where $u^{j}$ is the solution to $(\ref{1.1})$ that scatters backward in time to $S(t)(\phi_{0}^{j}, \phi_{1}^{j})$, and the remainder may be absorbed into $(R_{0, n}^{N}, R_{1, n}^{N})$. In this case as well,
\begin{equation}\label{6.20}
\| u^{j} \|_{L_{t,x}^{2(p - 1)}(\mathbf{R} \times \mathbf{R}^{3})} \leq M_{j} < \infty.
\end{equation}
The proof for $\lambda_{n}^{j} t_{n}^{j} \searrow -\infty$ is similar.

Also, by $(\ref{6.8})$, there are only finitely many $j$ such that $\| \phi_{0}^{j} \|_{\dot{H}^{s_{c}}} + \| \phi_{1}^{j} \|_{\dot{H}^{s_{c} - 1}} > \epsilon$. For all other $j$, small data arguments imply
\begin{equation}\label{6.21}
\| u^{j} \|_{L_{t,x}^{2(p - 1)}(\mathbf{R} \times \mathbf{R}^{3})} \lesssim \| \phi_{0}^{j} \|_{\dot{H}^{s_{c}}} + \| \phi_{1}^{j} \|_{\dot{H}^{s_{c} - 1}}.
\end{equation}
Then make use of the perturbative lemma.
\begin{lemma}[Perturbation lemma]\label{l6.2}
Let $I \subset \mathbf{R}$ be a time interval. Let $t_{0} \in I$, $(u_{0}, u_{1}) \in \dot{H}^{s_{c}} \times \dot{H}^{s_{c} - 1}$ and some constants $M$, $A$, $A' > 0$. Let $\tilde{u}$ solve the equation
\begin{equation}\label{6.22}
(\partial_{tt} - \Delta) \tilde{u} = F(\tilde{u}) = e,
\end{equation}
on $I \times \mathbf{R}^{3}$, and also suppose $\sup_{t \in I} \| (\tilde{u}(t), \partial_{t} \tilde{u}(t)) \|_{\dot{H}^{s_{c}} \times \dot{H}^{s_{c} - 1}} \leq A$, $\| \tilde{u} \|_{L_{t,x}^{2(p - 1)}(I \times \mathbf{R}^{3})} \leq M$,
\begin{equation}\label{6.23}
\| (u_{0} - \tilde{u}(t_{0}), u_{1} - \partial_{t} \tilde{u}(t_{0})) \|_{\dot{H}^{s_{c}} \times \dot{H}^{s_{c} - 1}} \leq A',
\end{equation}
and
\begin{equation}\label{6.24}
\| e \|_{L_{t}^{\frac{2}{1 + s_{c}}} L_{x}^{\frac{2}{2 - s_{c}}}(I \times \mathbf{R}^{3})} + \| S(t - t_{0})(u_{0} - \tilde{u}(t_{0}), u_{1} - \partial_{t} \tilde{u}(t_{0})) \|_{L_{t,x}^{2(p - 1)}(I \times \mathbf{R}^{3})} \leq \epsilon.
\end{equation}
Then there exists $\epsilon_{0}(M, A, A')$ such that if $0 < \epsilon < \epsilon_{0}$ then there exists a solution to $(\ref{1.1})$ on $I$ with $(u(t_{0}), \partial_{t} u(t_{0})) = (u_{0}, u_{1})$, $\| u \|_{L_{t,x}^{2(p - 1)}(I \times \mathbf{R}^{3})} \leq C(M, A, A')$, and for all $t \in I$,
\begin{equation}\label{6.25}
\| (u(t), \partial_{t} u(t)) - (\tilde{u}(t), \partial_{t} \tilde{u}(t)) \|_{\dot{H}^{s_{c}} \times \dot{H}^{s_{c} - 1}} \leq C(A, A', M)(A' + \epsilon).
\end{equation}
\end{lemma}
By Lemma $\ref{l6.2}$, the asymptotic orthogonality property $(\ref{6.7})$, and $(\ref{6.21})$,
\begin{equation}\label{6.26}
\limsup_{n \nearrow \infty} \| u^{n} \|_{L_{t,x}^{2(p - 1)}(\mathbf{R} \times \mathbf{R}^{3})}^{2} \lesssim \sum_{j} \| u^{j} \|_{L_{t,x}^{2(p - 1)}(\mathbf{R} \times \mathbf{R}^{3})}^{2} < \infty.
\end{equation}
This proves Theorem $\ref{t1.2}$. $\Box$

\bibliographystyle{plain}

\begin{thebibliography}{[00]}
\bibitem{BG}
	\newblock H. Bahouri and P. G\'erard,
	\newblock ``High frequency approximation of solutions to critical nonlinear wave equations."
	\newblock \textit{American Journal of Mathematics} \textbf{121} (1999), no. 1, 131 -- 175.

\bibitem{BS}
	\newblock H. Bahouri and J. Shatah,
	\newblock ``Decay estimates for the critical semilinear wave equation."
	\newblock \textit{Ann. Inst. H. Poincare Anal. Non Lineaire} \textbf{15} (1998), no. 6, 783--789.
	
\bibitem{D2}
	\newblock B. Dodson,
	\newblock ``Global well-posedness and scattering for the radial, defocusing, cubic wave equation with initial data in a critical Besov space",
	\newblock Preprint, arXiv:1608.02020. To appear, Analysis and PDEs.	
	
\bibitem{D}
	\newblock B. Dodson,
	\newblock ``Global well-posedness and scattering for the radial, defocusing, cubic nonlinear wave equation",
	\newblock Preprint,  arXiv:1809.08284.
	
\bibitem{GV}
	\newblock J. Ginibre and G. Velo
	\newblock ``Generalized {S}trichartz inequalities for the wave equation."
	\newblock \textit{Journal of Functional Analysis} \textbf{133} (1995), no. 1, 50 -- 68.		

\bibitem{GSV}
	\newblock J. Ginibre, A. Soffer, and G. Velo,
	\newblock ``The global Cauchy problem for the critical nonlinear wave equation."
	\newblock \textit{J. Funct. Anal.} \textbf{110} (1992), no. 1, 96--130.

\bibitem{Gril}
	\newblock M. Grillakis,
	\newblock ``Regularity and asymptotic behaviour of the wave equation with critical nonlinearity."
	\newblock \textit{Annals of Mathematics} \textbf{132} (1990), 485--509.
	
\bibitem{KlMa}
	\newblock S. Klainerman and M. Machedon,
	\newblock ``Space-times estimates for null forms and the local existence theorem."
	\newblock \textit{Comm. Pure Appl. Math.} \textbf{46} (1993), no. 9, 1221--1268.
	
\bibitem{KPV}
	\newblock C. Kenig, G. Ponce, and L. Vega,
	\newblock ``Global well-posedness for semi-linear wave equations",
	\newblock \textit{Comm. Partial Differential Equations} \textbf{25} (2000), no. 9--10, 1741--1752.	
	
\bibitem{Ramos}
	\newblock J. Ramos,
	\newblock ``A refinement of the Strichartz inequality for the wave equation with applications."
	\newblock \textit{Advances in Mathematics} \textbf{230} (2012), 649 -- 698.			
	
\bibitem{Shatah - Struwe}
	\newblock J. Shatah and M. Struwe,
	\newblock ``Regularity results for nonlinear wave equations."
	\newblock \textit{Ann. of Math. (2)} \textbf{138} (1993), no. 3, 503--518.

\bibitem{Shen}
	\newblock R. Shen,
	\newblock ``On the energy subcritical, nonlinear wave equation in $\mathbf{R}^{3}$ with radial data."
	\newblock \textit{Anal. PDE} \textbf{6} (2013), no. 8, 1929--1987.
	
\bibitem{Shen1}
	\newblock R. Shen,
	\newblock ``Scattering of solutions to the defocusing energy subcritical semi-linear wave equation in 3{D}."
	\newblock \textit{Comm. Partial Differential Equations} \textbf{42} (2017), no. 4, 495--518.

\bibitem{Sterb}
	\newblock J. Sterbenz,
	\newblock ``Angular regularity and Strichartz estimates for the wave equation," with an appendix by Igor Rodnianski.
	\newblock \textit{Int. Math. Res. Not.} (2005), no. 4, 187--231.	
	
\bibitem{Stri}
	\newblock R. S. Strichartz,
	\newblock ``Restrictions of {F}ourier transforms to quadratic surfaces and decay of solutions of wave equations.''
	\newblock \textit{Duke Mathematical Journal} \textbf{44} no. 3 (1977) 705 - 714.		

\bibitem{Struwe}
	\newblock M. Struwe,
	\newblock ``Globally regular solutions to the $u^{5}$ Klein - Gordon equation.''
	\newblock \textit{Ann. Scuola Norm. Sup. Pisa Cl. Sci. (4)} \textbf{15} (1988), no. 3, 495--513.
	
\bibitem{Tao}
	\newblock T. Tao,
	\newblock ``Spacetime bounds for the energy-critical nonlinear wave equation in three spatial dimensions."
	\newblock \textit{Dyn. Partial Differ. Equ.} \textbf{3} (2006), no. 2, 93--110.
	


\end{thebibliography}

\end{document}